\documentclass{article}

\title{Specialization of Appell's functions to univariate hypergeometric functions}

\author{Raimundas Vid\=unas\\ 
\em Kobe University\footnote{Address: Organization of Advanced Science and Technology, 
Kobe University, Rokko-dai 1-1, Nada-ku, Kobe 657-8501, Japan. 
Email: {\sf vidunas@math.kobe-u.ac.jp}. Tel: +81-78-803-5611. Fax: +81-78-803-5610.}}

\newtheorem{theorem}{Theorem}[section]

\newtheorem{definition}[theorem]{Defintion}

\usepackage{latexsym}

\newcommand{\app}[4]{F_{\!#1}\!
  \left(\left.{#2 \atop #3}\right| #4 \right) }

\newcommand{\hpg}[5]{{}_{#1}\mbox{\rm F}_{\!#2}\!
  \left(\left.{#3 \atop #4}\right| #5 \right) }

\newcommand{\hpgo}[2]{{}_{#1}\mbox{\rm F}_{\!#2}}

\newcommand{\proof}{{\bf Proof. }}
\newcommand{\qed}{\hfill $\Box$}
\newcommand{\equal}{&\!\!=\!\! &}

\newcommand{\CC}{\mbox{\bf C}}

\newcommand{\PP}{\mbox{\bf P}}
\newcommand{\ZZ}{\mbox{\bf Z}}

\date{}

\begin{document}

\maketitle

\begin{abstract}
Univariate specializations of Appell's hypergeometric functions $F_1$, $F_2$, $F_3$, $F_4$ satisfy ordinary Fuchsian equations of order at most 4. In special cases, these differential equations are of order 2 and could be simple (pullback) transformations of Euler's differential equation for the Gauss hypergeometric function. The paper classifies these cases, and presents corresponding relations between univariate specializations of Appell's functions and univariate hypergeometric functions. The computational aspect and interesting identities are discussed.
\end{abstract}

\section{Introduction}

Appell's hypergeometric functions $F_1(x,y)$, $F_2(x,y)$, $F_3(x,y)$, $F_4(x,y)$ are defined by the following double hypergeometric series: 
\begin{eqnarray} \label{appf1}
\app1{a;\;b_1,b_2}{c}{x,\,y} \equal \sum_{n=0}^{\infty} \sum_{m=0}^{\infty}
\frac{(a)_{n+m}\,(b_1)_n\,(b_2)_m}{(c)_{n+m}\;n!\,m!}\,x^n\,y^m,\\ \label{appf2}
\app2{a;\;b_1,b_2}{c_1,c_2}{x,\,y} \equal \sum_{n=0}^{\infty} \sum_{m=0}^{\infty}
\frac{(a)_{n+m}\,(b_1)_n\,(b_2)_m}{(c_1)_n\,(c_2)_m\;n!\,m!}\,x^n\,y^m,\\ \label{appf3}
\app3{\!a_1,a_2;\,b_1,b_2}{c}{x,\,y} \equal \sum_{n=0}^{\infty} \sum_{m=0}^{\infty}
\frac{(a_1)_n(a_2)_m(b_1)_n(b_2)_m}{(c)_{n+m}\;n!\,m!}\,x^n\,y^m,\\ \label{appf4}
\app4{a;\;b}{c_1,c_2\,}{x,\,y} \equal \sum_{n=0}^{\infty} \sum_{m=0}^{\infty}
\frac{(a)_{n+m}\,(b)_{n+m}}{(c_1)_n\,(c_2)_m\;n!\,m!}\,x^n\,y^m.
\end{eqnarray}
They are bivariate generalizations of the Gauss hypergeometric series
\begin{equation} \label{gausshpg}
\hpg21{A,\,B}{C}{z} = \sum_{n=0}^{\infty} 
\frac{(A)_{n}\,(B)_n}{(C)_n\,n!}\,z^n.
\end{equation}
In particular, Appell's functions satisfy Fuchsian  
(that is, linear homogeneous and with regular singularities) 
systems of partial differential equations, 
holonomic of rank 3 or 4, 
which are analogous to Euler's hypergeometric equation for the $\hpgo21$ function.
Euler's equation for (\ref{gausshpg}) is
\begin{equation} \label{eq:euler}
z(1-z)\,\frac{d^2y(z)}{dz^2}+
\big(C-(A+B+1)z\big)\frac{dy(z)}{dz}-A\,B\,y(z)=0.
\end{equation}
It is a Fuchsian equation with three singularities: $z=0$, $z=1$ and $z=\infty$. Local exponent differences are equal, respectively, to $1-C$, $C-A-B$, $A-B$ at them.

When the two arguments $x,y$ of Appell's functions are algebraically related, 
the univariate specializations satisfy Fuchsian ordinary
differential equations of order at most 4. This paper considers the following questions: 
\begin{quote}
Which univariate specializations of Appell's functions satisfy a second order Fuchsian 
ordinary differential equation? 

Which univariate specializations of Appell's functions satisfy Euler's equation (\ref{eq:euler})
for some parameters $A$, $B$, $C$, up to a projective transformation
\mbox{$y(z)\mapsto\theta(z) y(z)$} 
or a more general pull-back transformation 
\begin{equation}
z\mapsto\varphi(x), \qquad y(z)\mapsto\theta(x) y(\varphi(x)),
\end{equation}
where $\theta(z)$ or $\theta(x)$ is a power function, and $\varphi(x)$ is a rational function?
\end{quote}
Other interesting questions are: Which univariate specializations of Appell's functions satisfy
a Fuchsian ordinary differential equation of order 3? Up to pull-back transformations, which univariate specializations of Appell's functions satisfy differential equations for generalized hypergeometric functions $\hpgo32(z)$ or $\hpgo43(z)$?

The paper classifies univariate specializations of Appell's $F_2$, $F_3$, $F_4$ functions satisfying second order ordinary Fuchsian equations, when those ordinary equations follow from the partial differential equations for Appell's functions without reductions due to factorization of the respective differential operators (or reducibility of the monodromy group).
In particular, we do not consider isolated solutions of reducible Appell's systems of partial differential equations that are expressible via Gauss hypergeometric functions. Definition \ref{def:whole} below specifies the way ordinary differential equations {\em follow fully} from systems of partial differential equations.

We identify those univariate specialization of Appell's $F_2$, $F_3$, $F_4$ functions that can be expressed in terms of Gauss' hypergeometric function (\ref{gausshpg}). Their ordinary differential
equations are pull-back transformations of Euler's equation (\ref{eq:euler}), and {\em follow fully}
from the respective partial differential equations. 
To deal with univariate specializations onto singularity curves of the respective systems of partial differential equations, we present general ordinary Fuchsian equations for them. 
In the cases of $F_1$, $F_2$, $F_3$ functions, the specializations onto singularity curves can be expressed in terms of $\hpgo21$ or $\hpgo32$ hypergeometric functions. 

In Appendix, we briefly review identities between Appell's and Gauss hypergeometric functions
existing in literature. Our explicit results that appear more or less novel are the following:
\begin{itemize}
\item A case of $F_2(x,2-x)$ function that can be expressed as a $\hpgo21$ function with a quadratic argument; see Theorem \ref{th:f2cases}.  Equivalent identification 
of \mbox{$F_3(x,x/(2x-1))$} functions was noticed by Karlsson  \cite{Karlsson80}.
These relations apply nicely when the monodromy group of the $\hpgo21$ function is dihedral;
then the $F_2$ and $F_3$ solutions of the same ordinary Fuchsian equation are terminating,
and dihedral $\hpgo21$ functions can be expressed as elementary functions. 
We show this application in Section \ref{sec:dihedral}.
\item A separation of variables case for the $F_2$ function; see formula (\ref{eq:f2separation}).
It is related to the well-known Bailey case \cite{Bailey33} of variable separation
for the $F_4$ function
\begin{equation} \label{eq:bailey}
\app4{a;\;b}{c,a\!+\!b\!-\!c\!+\!1}{x(1-y),y(1-x)}
=\hpg21{a,b}{c}{x} \hpg21{a,\;b}{a\!+\!b\!-\!c\!+\!1}{y}
\end{equation}
via a known transformation between $F_2$ and $F_4$ functions;
\item Translation of the above two $F_2$ cases to relations of $F_1$ and $F_3$ functions
to Gauss hypergeometric functions.
\item A few cases of $F_4(t^2,(1-t)^2)$ functions (that is, specializations to the quadratic singular curve
of $F_4$) expressible as $\hpgo21$ or $\hpgo32$ functions, or products of two $\hpgo21$ solutions
of the same Euler's equation (\ref{eq:euler}). 
\end{itemize}

Identities between bivariate and univariate hypergeometric series are usually derived by methods of series manipulation or using integral representations. The method of relating hypergeometric functions as solutions of coinciding differential equations is usually considered tedious and computationally costly. With powerful computer algebra techniques available, the method of identifying differential equations can be worked out rather comprehensively.
We discuss general computational techniques in Section \ref{sec:compute}.  Here below and in Section \ref{sec:f2} we gradually introduce the algebraic setting, computational methods and some shortcuts.

Consider a holomorphic function $F(x,y)$ on an open set in $\CC\times\CC$, 
and a univariate specialization $F(x(t),y(t))$ of it. The full derivatives with respect to $t$ 
are expressed linearly in terms of the partial derivatives $\partial F/\partial x$, $\partial F/\partial y$, 
$\partial^2F/\partial x^2$, etc. For example,
\begin{eqnarray} \label{eq:f2d1g}
\frac{dF}{dt} \equal \dot{x}\,\frac{\partial F}{\partial x}+\dot{y}\,\frac{\partial F}{\partial y},\\
\label{eq:f2d2g}\frac{d^2F}{dt^2}\! \equal\! \dot{x}^{2}\,\frac{\partial^2 F}{\partial x^2}+
2\dot{x}\dot{y}\,\frac{\partial^2 F}{\partial x\partial y}+
\dot{y}^{2}\,\frac{\partial^2 F}{\partial y^2}
+\ddot{x}\,\frac{\partial F}{\partial x}+\ddot{y}\,\frac{\partial F}{\partial y}.
\end{eqnarray}
Here $\dot{x}=dx/dt$, $\dot{y}=dy/dt$, $\ddot{x}=d^2x/dt^2$, $\ddot{y}=d^2y/dt^2$.
To compute next order full derivatives, one applies the Leibniz rule and lets $d/dt$ act
on the partial derivatives by copying the action on $F$ in (\ref{eq:f2d1g}). 

We identify partial differential equations with ordinary differential equations in 
the specialization variable $t$ by identifying their ``mixed" forms in the partial derivatives with
the coefficients specialized to functions in $t$.
To fix a univariate specialization mapping, consider a holomorphic map $\Phi$ from an open subset
of $\CC$ to $\CC\times\CC$, mapping $t\mapsto(x(t),y(t))$. If $F(x,y)$ is a holomorphic function
on the image of $\Phi$, then $F\circ\Phi$ is the univariate specialization $F(x(t),y(t))$.
Geometrically, a univariate specialization is the pullback of $F(x,y)$ with respect to $\Phi$.

\begin{definition} \label{def:whole} \rm
If $E$ is an ordinary differential equation with respect to $t$, 
its {\em partial differential form} under $\Phi$ 
is the expression where the derivatives with respect to $t$ are replaced
by partial derivatives following formulas (\ref{eq:f2d1g})--(\ref{eq:f2d2g}), etc.

If $\widetilde{E}$ is a partial differential equation with respect to $x$, $y$, 
its {\em specialized form} under $\Phi$ 
is the expression where the coefficients to the partial derivatives
are specialized $x\mapsto x(t)$, $y\mapsto y(t)$.

Let $H$ denote a system of partial differential equations with respect to $x$, $y$.
The ordinary differential equation $E$ is said to {\em follow fully} from $H$ under $\Phi$, 
if its partial differential form under $\Phi$ coincides with the specialized form
(under $\Phi$) of some partial differential equation following from $H$ by algebraic
and partial differentiation operations.

We also say that $E$ is {\em implied fully} by $H$, if there is a specialization map $\Phi$ 
such that $E$ follows fully from $H$ under $\Phi$. 
Here we allow $H$, $E$ and $\Phi$ to have parameters; the parameters of $H$ may specialize
(to constants or functions in the parameters of $E$, $\Phi$) in the specialized forms 
of its partial differential equations under $\Phi$.
\end{definition}

The classical systems of partial differential equations for Appell's functions
$F_2$, $F_3$, $F_4$, $F_1$ are presented, respectively, in formulas (\ref{app2a})--(\ref{app2b}), 
(\ref{app3a})--(\ref{app3b}), (\ref{app4a})--(\ref{app4b}), \mbox{(\ref{app1a})--(\ref{app1b})} below.
Algebraically, linear partial differential equations with polynomial coefficients are identified
with partial differential operators 
in the Weyl algebra \mbox{$\CC[x,y]\langle\partial/\partial x,\partial/\partial y\rangle$}.
A system of partial differential equations corresponds to a left ideal in the Weyl algebra.

The systems of partial differential equations for $F_2$, $F_3$, $F_4$ functions have rank $4$.
Correspondingly, Gr\"obner bases for the mentioned left ideals in the Weyl algebra 
give linear expressions of all partial derivatives in terms of 4 partial derivatives of low order,
say $\partial^2F/\partial x^2$ 
$\partial F/\partial x$, $\partial F/\partial y$, $F$.
If we take a total degree ordering of partial derivatives, the leading coefficients in
the expressions for higher order derivatives vanish only when specialized
onto a singularity curve of the differential system.

If a univariate specialization of Appell's $F_2$, $F_3$ or $F_4$ function is not onto
a singularity curve of the differential system,
the full derivatives can be expressed linearly in terms of the 4 basic partial derivatives,
and a linear relation between $d^4F/dt^4$,  $d^3F/dt^3$, $d^2F/dt^2$, $dF/dt$, $F$ 
gives an ordinary differential equation of order at most 4 for the univariate function. 
If the univariate specialization is onto a singularity curve, it turns out that we get
an ordinary differential equation of order less than $4$. 
All these ordinary differential equations are Fuchsian.

Similarly, the system of partial differential equations for the $F_1$ function has rank $3$, 
and its univariate specializations satisfy ordinary Fuchsian equations of order at most $3$. 

We are interested in cases when the ordinary differential equation has order 2. 
Such a second order differential equation must follow from partial differential equations
of order at most 2. For non-singular specializations of $F_2$, $F_3$, $F_4$ functions,
linear relations between partial derivatives of order at most 2 are generated by the classical pairs
of partial differential equations. 

\section{Identities with Appell's $F_2$ function}
\label{sec:f2}

A system partial differential equations for the $F_2(x,y)$ function is:
\begin{eqnarray} \label{app2a}
x(1-x)\frac{\partial^2F}{\partial x^2}-xy\frac{\partial^2F}{\partial x\partial y}
+\left(c_1-(a+b_1+1)x\right)\frac{\partial F}{\partial x}-b_1y\frac{\partial F}{\partial y}
-ab_1F=0,\\ \label{app2b}
y(1-y)\frac{\partial^2F}{\partial y^2}-xy\frac{\partial^2F}{\partial x\partial y}
+\left(c_2-(a+b_2+1)y\right)\frac{\partial F}{\partial y}-b_2x\frac{\partial F}{\partial x}
-ab_2F=0.
\end{eqnarray}
This is an equation system of rank 4. 
Its singular locus on $\PP^1\times\PP^1$ is the union of the following lines:
\begin{equation} \label{app2sing}
x=0,\quad x=1,\quad x=\infty, \quad y=0,\quad y=1,\quad y=\infty, \quad x+y=1.
\end{equation}
First we handle univariate specializations onto the singular lines. We determine that
generally these specializations satisfy hypergeometric equations for $\hpgo21$ or $\hpgo32$
functions. Explicit identification of the following three function pairs is presented after
Theorem \ref{th:f2cases}. Those identities can be easily derived from power series
manipulations as well.
\begin{theorem} \label{th:hpg32}
The following function pairs satisfy the same ordinary differential 
equations (of order $2$ or $3$):
\begin{itemize}
\item $\displaystyle\app2{a;\;b_1,b_2}{c_1,c_2}{x,\,0}$ and $\displaystyle\hpg21{a,\,b_1}{c_1}{x}$;
\item $\displaystyle\app2{a;\;b_1,b_2}{c_1,c_2}{x,\,1}$ and 
$\displaystyle\hpg32{a,\,b_1,\,a-c_2+1}{c_1,\,a+b_2-c_2+1}{x}$;\\
\item $\displaystyle\app2{a;\;b_1,b_2}{c_1,c_2}{x,1-x}$ and
$\displaystyle(1-x)^{-a}\,\hpg32{a,\,c_1-b_1,\,a-c_2+1}{c_1,\,a+b_2-c_2+1}{\frac{x}{x-1}}$.
\end{itemize}
\end{theorem} 
\proof The first specialization is trivial: after the substitution $y=0$ into the coefficients
of (\ref{app2a}) we have differentiation with respect to $x$ only.

To handle the other two cases, we use the third order differential equation for the
$\hpgo32(x)$ function; see formula (\ref{eq:hpgde3}) in the Appendix.

For the $F_2(x,1)$ function, subtract (\ref{app2b}) from (\ref{app2a}) and set $y=1$
in the coefficients:
\begin{eqnarray} \label{eq:dy1}
x(1-x)\frac{\partial^2F}{\partial x^2}
+\left(c_1-(a+b_1-b_2+1)x\right)\frac{\partial F}{\partial x} \nonumber\\
-\left(c_2-a+b_1-b_2-1\right)\frac{\partial F}{\partial y}
-a(b_1-b_2)F&=&0.
\end{eqnarray}
Then make the following combination of partial differential equations for the $F_2(x,1)$: 
\begin{equation}
x\frac{\partial}{\partial x}\mbox{[Eq. (\ref{eq:dy1})]}+(a+b_2-c_2+1)\mbox{[Eq. (\ref{app2a})]}
-b_1\mbox{[Eq. (\ref{app2b})]}.
\end{equation}
After setting $y=1$ in the coefficients, we get an equation with differentiation by 
$x$ only. That differential equation is identified as hypergeometric equation
(\ref{eq:hpgde3}) with the indicated values of the parameters $A,B,C,D,E$.

For the $F_2(x,1-x)$ function, add (\ref{app2a}) and (\ref{app2b}) and set $y=1-x$:
\begin{eqnarray} \label{eq:sum2}
x(1-x)\left( \frac{\partial^2F}{\partial x^2}-2\frac{\partial^2F}{\partial x\partial y}
+\frac{\partial^2F}{\partial y^2} \right)+
\left(c_1-(a+b_1+b_2+1)x\right)\frac{\partial F}{\partial x} \nonumber\\
+\left(c_2-a-b_1-b_2-1+(a+b_1+b_2+1)x\right)\frac{\partial F}{\partial y}
-a(b_1+b_2)F&=&0,
\end{eqnarray}
and consider the following combination of partial differential equations for it:
\begin{equation} \label{eq:comb2sum1}
\left(\frac{\partial}{\partial x}-\frac{\partial}{\partial y}+\frac{a+1}{x-1}\right)\mbox{[Eq. (\ref{eq:sum2})]}
+\frac{a+b_2-c_2+1}{x\,(1-x)}\mbox{[Eq. (\ref{app2a})]}
+\frac{c_1-b_1}{x(1-x)}\mbox{[Eq. (\ref{app2b})]}.
\end{equation}
With $y=1-x$, we recognize the full derivatives
\begin{eqnarray} 
&& \frac{d^3F}{dx^3}=\frac{\partial^3F}{\partial x^3}
-3\frac{\partial^3F}{\partial x^2\partial y}+3\frac{\partial^3F}{\partial x\partial y^2}
-\frac{\partial^3F}{\partial y^3},\\ \label{partial2}
&& \frac{d^2F}{dx^2}=\frac{\partial^2F}{\partial x^2}-2\frac{\partial^2F}{\partial x\partial y}
+\frac{\partial^2F}{\partial y^2},\qquad
\frac{dF}{dx}=\frac{\partial F}{\partial x}-\frac{\partial F}{\partial y}
\end{eqnarray}
in (\ref{eq:comb2sum1}), and get  an ordinary differential equation of order 3, 
with the following singularities and  local exponents:
\begin{eqnarray*}
\begin{array}{llcl}
\mbox{at $x=0$}: & 0, & 1-c_1, & c_2-a-b_2;\\
\mbox{at $x=1$}: & 0, & 1-c_2, & c_1-a-b_1;\\
\mbox{at $x=\infty$}: & a, & a+1, & b_1+b_2.
\end{array}
\end{eqnarray*}
The differential equation for the $\hpgo32$ function (times the power factor) has the same data,
making the requisite check of the differential equation on the second function worthwhile.
\qed\\

The cases when $\hpgo32$ functions in Theorem \ref{th:hpg32} become $\hpgo21$ functions
can be found by equating a pair of upper and lower parameters. 
Note that the case $c_2=a-b_1+b_2+1$ for the $F_2(x,1)$ function is already evident from (\ref{eq:dy1}),
while the case $c_1+c_2=a+b_1+b_2+1$ for the $F_2(x,1-x)$ function is visible in (\ref{eq:sum2}).
Euler's equation (\ref{eq:euler}) is easily recognizable in those two formulas.

The next two theorems tell that second order Fuchsian equations for univariate $F_2$
specializations outside the singularity lines, implied fully by \rm(\ref{app2a})--(\ref{app2b}),
follow from \mbox{\rm(\ref{app2a})--(\ref{app2b})} without differentiation. 
\begin{theorem} \label{th:app2o2}
The following univariate specializations of Appell's $F_2(x,y)$ function satisfy
ordinary Fuchsian equations of second order (with respect to $x$ or $t$):
\begin{eqnarray} \label{eq:f2f0}
\app2{a;\;b_1,b_2}{2b_1,2b_2}{x,2-x},\qquad
\label{eq:f2f5} \app2{b_1+b_2-\frac12;\;b_1,b_2}{2b_1,2b_2}{1-t^2,1-\frac{(t+s)^2}{s^2-1}}.
\end{eqnarray}
In the second function, $s$ is any constant.
\end{theorem}
\proof For the $F_2(x,2-x)$ function, divide (\ref{app2a}) by $x$ and (\ref{app2b}) by $-y$, 
add the two equations, and substitute $y=2-x$ in the coefficients:
\begin{eqnarray}
(1-x)\left( \frac{\partial^2F}{\partial x^2}-2\frac{\partial^2F}{\partial x\partial y}
+\frac{\partial^2F}{\partial y^2} \right)
+\left( \frac{c_1}{x}-a-b_1-1-\frac{b_2x}{x-2} \right)\frac{\partial F}{\partial x} \nonumber\\
+\left( \frac{c_2}{x-2}+a+b_2+1-b_1\frac{2-x}{x} \right)\frac{\partial F}{\partial y}
-a\left( \frac{b_1}{x}+\frac{b_2}{x-2} \right)F&=&0.
\end{eqnarray}
If $c_1=2b_1$ and $c_2=2b_2$, following (\ref{partial2}) we recognize the second order equation 
\begin{equation} \label{eq:xy2}
(1-x)\frac{d^2F}{dx^2}
+\left( \frac{2b_1}{x}-\frac{2b_2}{x-2}-a-b_1-b_2-1 \right)\frac{dF}{dx}
-a\left( \frac{b_1}{x}+\frac{b_2}{x-2} \right)F=0.
\end{equation}
For the second $F_2$ function, the following differential equation can be checked
by substituting (\ref{eq:f2d1g})--(\ref{eq:f2d2g}) with evaluated $\dot{x}$, $\dot{y}$, 
$\ddot{x}$, $\ddot{y}$, then linearly eliminating the derivatives 
$\partial^2 F/\partial x^2$, $\partial^2 F/\partial y^2$ using equations (\ref{app2a})--(\ref{app2b}):
\begin{eqnarray} \label{eq:f2sep}
\frac{d^2F}{dt^2}
+\!\left(\frac{4b_1t}{t^2-1}+\frac{4b_2(t+s)}{t^2+2st+1} \right)\!\frac{dF}{dt}
+\left({\textstyle b_1+b_2-\frac12}\right)\!\left(\frac{4b_1}{t^2-1}+\frac{4b_2}{t^2+2st+1}\right)\!F =0.
\end{eqnarray}
\qed

\begin{theorem}
Suppose that a univariate specialization of Appell's $F_2(x,y)$ function satisfies a second order 
ordinary Fuchsian equation fully implied by the system \mbox{\rm(\ref{app2a})--(\ref{app2b})}.
Then either the specialization is into singular locus \mbox{\rm(\ref{app2sing})}, 
or they are represented by one of the two functions in \mbox{\rm(\ref{eq:f2f0})}.
\end{theorem}
\proof  We can assume that $x(x-1)y(y-1)(x+y-1)\neq 0$.
If $F(t)=F_2\left(x(t),y(t)\right)$ satisfies a second order Fuchsian equation,
second order partial derivatives come from the $d^2F/dt^2$ term following (\ref{eq:f2d2g}).
The rank of the differential system (\ref{app2a})--(\ref{app2b}) is 4, and 
there are 2 linearly independent relations between partial derivatives of order at most 2.
Therefore equations (\ref{app2a})--(\ref{app2b}) linearly generate all partial differential 
equations of order 2  that algebraically follow from them. 
(In other words, the corresponding two differential operators in 
the Weyl algebra \mbox{$\CC[x,y]\!<\!\!\partial/\partial x,\partial/\partial y\!\!>$} 
linearly generate order 2 operators of the corresponding left ideal.)

Linear elimination of $\partial^2 F/\partial x^2$, $\partial^2 F/\partial y^2$ from (\ref{eq:f2d2g})
using equations  (\ref{app2a})--(\ref{app2b}) gives the following coefficient to 
$\partial^2 F/\partial x\partial y$:
\begin{equation} \label{eq:genf2}
\frac{y}{1-x}\,\dot{x}^{2}+2\dot{x}\dot{y}
+\frac{x}{1-y}\,\dot{y}^{2}.
\end{equation}
This expression must be equated to $0$. 
The expression is a quadratic form in the derivatives $\dot{x},\dot{y}$, 
with the discriminant equal to $4(1-x-y)/(1-x)(1-y)$. We can factorize (\ref{eq:genf2}) into linear 
differential forms if we parametrize the surface $(1-x-y)u^2=(1-x)(1-y)$. 
Accordingly, we substitute
\begin{equation} \label{eq:solvy}
y=\frac{(1-u^2)(x-1)}{x+u^2-1}.
\end{equation}
The parameter $u$ is undefined only if $x+y=1$.

Expression (\ref{eq:genf2}) then factorizes as follows:
\begin{equation}  \label{eq:facteq}
\frac{\big(2x(x-1)du-(u-1)(x-u-1)dx\big)\big(2x(x-1)du-(u+1)(x+u-1)dx\big)}{\left(x+u^2-1\right)^3},
\end{equation}
where we write $du$, $dx$ instead of $du/dt$, $dx/dt$, because the variable $t$ is irrelevant. 
(It can be composed with any univariate function without changing our differential relations.)
The difference in the two numerator factors is in the sign of $u$.  
The differential equation
\begin{equation}
2x(x-1)du=(u-1)(x-u-1)dx
\end{equation}
is generally solved by the following algebraic relation between $u$ and $x$:
\begin{equation} \label{eq:intcurve}
(u+x-1)^2+C(x-1)(u-1)^2=0,
\end{equation}
where $C$ is an integration constant. After eliminating $u$ from (\ref{eq:solvy}) and  (\ref{eq:intcurve}) 
we get the algebraic relation between $x$ and $y$ that solves (\ref{eq:genf2}):
\begin{equation} \label{eq:intcurve2}
(x+y-2)^2-2C(y^2+xy-2x-2y+2)+C^2y^2=0.
\end{equation}
The algebraic curves (\ref{eq:intcurve}) and (\ref{eq:intcurve2}) can be parametrized as follows:
\begin{equation} \label{eq:f25para}
x=1-Ct^2, \qquad y=1-\frac{C(t+1)^2}{C-1}, \qquad u=\frac{Ct(t+1)}{Ct+1}.
\end{equation}
The parameter $t=(x+y-2)/2C-y/2$ is undefined if $C=0$ and $x+y=2$, which is a special case
of (\ref{eq:intcurve2}). We keep in mind the special case $x+y=2$. 

Now we have to check when the elimination of $\partial^2 F/\partial x^2$, $\partial^2 F/\partial y^2$
from (\ref{eq:f2d2g}) using equations (\ref{app2a})--(\ref{app2b}) and parametrization (\ref{eq:f25para})
leaves the quotient of the coefficients to $\partial F/\partial x$, $\partial F/\partial y$ equal to
$\dot{y}/\dot{x}$, so that we could get rid of the $\partial F/\partial x$, $\partial F/\partial y$ terms by adding a rational multiple of $dF/dt$ as in (\ref{eq:f2d1g}).
If $x+y=2$,  we are led to compare coefficients to first order partial derivatives in (\ref{eq:xy2}).
There remains the case (\ref{eq:f25para}) with $C\neq 0$.
We replace $C\mapsto s^2$, $t\mapsto t/s$ and get the parametrization
of $x$ and $y$ as in (\ref{eq:f2f5}). The derivative $d^2F/dx^2$ in (\ref{eq:f2d2g}) is reduced 
by (\ref{app2a})--(\ref{app2b}) to the following expression:
\begin{eqnarray*}
\left(\frac{2c_1t}{t^2-1}+\frac{4b_2(t+s)}{t^2+2st+1}+\frac{2a\!+\!2b_1\!-\!2b_2\!-\!2c_1\!+\!1}{t}\right)
\!2t\,\frac{\partial F}{\partial x}+\nonumber \hspace{108pt}\\
\left(\frac{4b_1t}{t^2\!-\!1}+\frac{2c_2(t+s)}{t^2\!+\!2st\!+\!1}
+\frac{2a\!-\!2b_1\!+\!2b_2\!-\!2c_2\!+\!1}{t+s}\right) 
\!\frac{2(t\!+\!s)}{s^2\!-\!1}\frac{\partial F}{\partial y} 
-\left(\frac{4ab_1}{t^2\!-\!1}+\frac{4ab_2}{t^2\!+\!2st\!+\!1}\right)\!F.
\end{eqnarray*}
To have the right quotient of coefficients to the first order partial derivatives,
we must have $c_1=2b_1$, $c_2=2b_2$ and $2a+1=2b_1+2b_2$.\qed\\

Now we indicate univariate hypergeometric solutions of the ordinary differential equations
implied by Theorem \ref{th:app2o2}.
\begin{theorem} \label{th:f2cases}
The following function pairs satisfy the same second order ordinary differential 
equations (with respect to $x$ or $t$):
\begin{itemize}
\item $\displaystyle\app2{a;\;b_1,b_2}{2b_1,2b_2}{x,2-x}$ and
$\displaystyle(x-2)^{-a}\,\hpg21{\frac{a}2,\frac{a+1}2-b_2}{b_1+\frac12}{\frac{x^2}{(2-x)^2}}$;\\
\item $\displaystyle\app2{b_1+b_2-\frac12;\;b_1,b_2}{2b_1,\,2b_2}
{-\frac{4t}{(t-1)^2},\frac{(1-s)(st^2-1)}{s(t-1)^2}}$ and\\
$\displaystyle(1-t)^{2b_1+2b_2-1}\hpg21{b_1+b_2-\frac12,b_2}{b_1+\frac12}{st^2}$. 
\end{itemize}
In the last pair, $s$ is any constant. 
\end{theorem}
\proof The differential equation (\ref{eq:xy2}) for the $F_2(x,2-x)$ function
has the following singularities and local exponents:
\begin{eqnarray*}
\begin{array}{lrr}
\mbox{at} \quad x=0: & 0, & 1-2b_1,\\
\mbox{at} \quad x=2: & 0, & 1-2b_2,\\
\mbox{at} \quad x=1: & 0, & b_1+b_2-a,\\
\mbox{at} \quad x=\infty: & a, & b_1+b_2.\end{array}
\end{eqnarray*}
Local exponent differences at $x=1$ and $x=\infty$ are both equal to $b_1+b_2-a$.
There is a chance that (\ref{eq:xy2}) is a pull-back of Euler's equation (\ref{eq:euler}) 
with respect to the covering $z\mapsto x^2/(2-x)^2$, with its local exponent differences
possibly the following:
\[\textstyle
1-C=\frac12-b_1, \qquad A-B=\frac12-b_2, 
\qquad C-A-B=b_1+b_2-a.
\]
The shift of the local exponents at $x=\infty$ and $x=2$ by $-a$ would be needed.
Therefore we consider Euler's equation (\ref{eq:euler}) 
with $C=1/2+b_1$, $A=a/2+1/2-b_2$, \mbox{$B=a/2$},
compute its pullback with respect to
the covering $z\mapsto x^2/(2-x)^2$ and projective normalization $y(x)\mapsto (x-2)^{-a}y(x)$,
and check that the obtained differential equation indeed coincides with (\ref{eq:xy2}).

To get the result for the second $F_2$ function, we first consider differential equation (\ref{eq:f2sep})
for the function (\ref{eq:f2f5}) of Theorem \ref{th:app2o2}. 
The differential equation has the following singularities and local exponents:
\begin{eqnarray*}
\begin{array}{lrr}
\mbox{at $t=1$ and $t=-1$}: & 0, & 1-2b_1,\\
\mbox{at the roots of $t^2+2st+1$}: 
 & 0, & 1-2b_2,\\
\mbox{at } t=\infty: & 2b_1+2b_2, & 2b_1+2b_2-1.
\end{array}
\end{eqnarray*}
If the point $t=\infty$ is an apparent singularity, equation (\ref{eq:f2sep}) might be a pull-back 
of Euler's equation (\ref{eq:euler}) with respect to the covering $z\mapsto K(t+1)^2/(t-1)^2$,
with the local exponent differences
\[\textstyle
1-C=\frac12-b_1, \qquad A-B=\frac12-b_1, 
\qquad C-A-B=1-2b_2,
\]
and the constant $K=(s+1)/(s-1)$ adjusted so that the other 2 singular points lie above $z=1$.
The shift of the local exponents at $t=\infty$ and $t=1$ by $1-2b_1-2b_2$ would be needed.

Indeed, this is a pullback. 
We conclude that the function in (\ref{eq:f2f5}) and 
\begin{equation}
(1-t)^{1-2b_1-2b_2}\hpg21{b_1+b_2-\frac12,b_2}{b_1+\frac12}
{\frac{(s+1)(t+1)^2}{(s-1)(t-1)^2}}
\end{equation}
satisfy the same second order differential equation. After substituting
$t\mapsto (t+1)/(t-1)$ and $s\mapsto (s+1)/(s-1)$ into the two functions, 
we obtain the last pair of functions of this theorem.
\qed\\

Recall \cite[Section 2.9]{Erdelyi81} that Euler's equation (\ref{eq:euler}) has 24 hypergeometric
$\hpgo21$ series solutions (representing 6 different functions) in general, as discovered by Kummer. 
Generalized hypergeometric equation (\ref{eq:hpgde3}) has 6 hypergeometric $\hpgo32$ 
solutions in general. Therefore one can take alternative univariate hypergeometric solutions
in Theorems \ref{th:hpg32} and \ref{th:f2cases}; the presented ones look most convenient
representatives. In Appendix, we recall similar sets of Appell's solutions, 
say (\ref{eq:f2s1})--(\ref{eq:f2s7}), to the classical systems  like (\ref{app2a})--(\ref{app2b}).

Now we discuss explicit relations between the function pairs in Theorems \ref{th:hpg32} and
\ref{th:f2cases}. In each case of Theorem \ref{th:hpg32}, the two 
functions differ by a constant multiple in a neighborhood of $x=0$ in general, 
because they satisfy the same Fuchsian differential equation, 
start with the same local power exponent at $x=0$, 
there are no other integer local exponents at $x=0$ in general, 
and both functions are defined by proper power series around $x=0$. 
In the first case, the two functions are equal to each other, trivially.
For the second function pair, it is tempting to expand the $F_2(x,1)$ function
using the power series definition (\ref{appf2}),
and use Gauss' formula \cite[Theorem 2.2.2]{AAR99} to evaluate the $\hpgo21(1)$ 
coefficients to the powers of $x$.
That would give the following formula:
\begin{eqnarray} \label{eq:wrongformula}
\app2{a;\;b_1,b_2}{c_1,c_2}{x,\,1} \mbox{``=" }
\frac{\Gamma(c_2)\Gamma(c_2-a-b_2)}{\Gamma(c_2-a)\Gamma(c_2-b_2)}
\hpg32{a,\,b_1,\,a-c_2+1}{c_1,\,a+b_2-c_2+1}{x},
\end{eqnarray}
Unfortunately, this identification is incorrect, because the $\hpgo21(1)$ coefficients
to high enough powers of $x$ will eventually diverge\footnote{In the published version of this paper,
in {\em J.~Math.~Anal.~Appl.} {\bf 355} (2009), pg.~145--163,
several wrong formulas of this type were claimed. The subsequent paper \cite{Vidunas09}
rectifies and explains the error.}; see \cite{Vidunas09}.

In the first case of Theorem \ref{th:f2cases}, the $F_2(x,2-x)$ series does not converge. 
Linear relations between $F_2(x,2-x)$ and $\hpgo21$ solutions of the same second order
Fuchsian equation are cumbersome in general, since analytic continuation has to be identified scrupulously.  

Nevertheless, if the parameters $b_1$ and $b_2$ in the \mbox{$F_2(x,2-x)$} series
are negative integers, then the series becomes a finite sum of hypergeometric terms, 
while $\hpgo21(x)$ functions have a dihedral monodromy group. 
The dihedral $\hpgo21(x)$ functions can be expressed elementarily in terms of nested radical functions, and the terminating $F_2$ series occur in the most general form of these elementary expressions. 
We discuss this 
in Section \ref{sec:dihedral}.

In the second case of Theorem \ref{th:f2cases}, let us substitute $s\mapsto s^2$, $t\mapsto t/s$. 
We conclude that 
\begin{equation} \label{eq:f2last}
\app2{b_1+b_2-\frac12;\;b_1,b_2}{2b_1,2b_2}
{-\frac{4st}{(t-s)^2},-\frac{(s^2-1)(t^2-1)}{(t-s)^2}}
\end{equation}
and
\begin{equation} \label{eq:f2hpgs}
(s-t)^{2b_1+2b_2-1}\hpg21{b_1+b_2-\frac12,b_2}{b_1+\frac12}{s^2}
\hpg21{b_1+b_2-\frac12,b_2}{b_1+\frac12}{t^2}
\end{equation}
satisfy the same second order ordinary differential equation, with respect to $t$.
In this setting, the $\hpgo21(s^2)$ expression is just a constant factor. 
But the symmetry between $s$ and $t$ suggests that equation system (\ref{app2a})--(\ref{app2b})
for $\app2{b_1+b_2-1/2;b_1,b_2}{2b_1,2b_2}{x,y}$ can be transformed following (\ref{eq:f2last}) 
to a differential system where the new variables $s$ and $t$ are separated.
Moreover, the two separated equations with respect to $s$ or $t$ 
would be simple transformations of Euler's equation (\ref{eq:euler}).
In Section \ref{sec:appellf4} we relate this situation by a known transformation
to the well-known case (\ref{eq:bailey}) of separation of variables for the $F_4(x,y)$ function.

Here is an identity illustrating the last case of Theorem \ref{th:f2cases}, obtainable after substitution
$s\mapsto (s+1)/(s-1)$ into (\ref{eq:f2last})--(\ref{eq:f2hpgs}) and consideration of $F_2$ and $\hpgo21$ solutions of the same differential system around the point $(t,s)=(0,0)$:
\begin{eqnarray} \label{eq:f2separation}
&&\app2{b_1+b_2-\frac12;\;b_1,b_2}{2b_1,2b_2}
{\frac{4(1-s^2)\,t}{(1+s+t-st)^2},\frac{4(1-t^2)\,s}{(1+s+t-st)^2}}= \nonumber\\
&&\left(\frac{1\!+\!s\!+\!t\!-\!st}{1-s}\right)^{2b_1+2b_2-1}
\hpg21{b_1\!+\!b_2\!-\!\frac12,b_2}{2b_2}{\frac{4s}{s^2\!-\!1}}
\hpg21{b_1\!+\!b_2\!-\!\frac12,b_2}{b_1+\frac12}{t^2}. \qquad%
\end{eqnarray}
One can check this identity by comparing low degree terms of
bivariate Taylor series around $(t,s)=(0,0)$.


\section{Identities with Appell's $F_3$ function}

A system partial differential equations for Appell's $F_3$ function is:
\begin{eqnarray} \label{app3a}
x(1-x)\frac{\partial^2F}{\partial x^2}+y\frac{\partial^2F}{\partial x\partial y}
+\left(c-(a_1+b_1+1)x\right)\frac{\partial F}{\partial x}-a_1b_1F=0,\\ \label{app3b}
y(1-y)\frac{\partial^2F}{\partial y^2}+x\frac{\partial^2F}{\partial x\partial y}
+\left(c_2-(a_2+b_2+1)y\right)\frac{\partial F}{\partial y}-a_2b_2F=0.
\end{eqnarray}
This is an equation system of rank 4.
The singular locus of this equation on $\PP^1\times\PP^1$ is the union of the following curves:
\begin{equation}
x=0,\quad x=1,\quad x=\infty, \quad y=0,\quad y=1,\quad y=\infty, \quad xy=x+y.
\end{equation}

There is a straightforward relation between $F_2(x,y)$ and $F_3(x,y)$ functions,
coming from a transformation between differential systems (\ref{app2a})--(\ref{app2b})
and (\ref{app3a})--(\ref{app3b}). The relation is that the functions
\begin{eqnarray} \label{app2app3}
\app2{a;\;b_1,b_2}{c_1,c_2}{x,y}, \qquad
x^{-b_1}\,y^{-b_2}\,\app3{1+b_1-c_1,1+b_2-c_2;\,b_1,b_2}{1+b_1+b_2-a}{\frac1x,\,\frac1y\,}.
\end{eqnarray}
both satisfy (\ref{app2a})--(\ref{app2b}). Therefore we can translate the results of 
Theorems \ref{th:hpg32} through \ref{th:f2cases} straightforwardly. 
\begin{theorem} \label{th:f3cases}
The following function pairs satisfy the same ordinary differential 
equations (or order $2$ or $3$, with respect to $x$ or $t$):
\begin{itemize}
\item $\displaystyle\app3{a_1,a_2;\,b_1,b_2}{c}{x,\,0}$ and $\displaystyle\hpg21{a_1,\,b_1}{c}{x}$;
\item $\displaystyle\app3{a_1,a_2;\;b_1,b_2}{c}{x,\,1}$ and
$\displaystyle\hpg32{a_1,b_1,c-a_2-b_2}{c-a_2,c-b_2}{x}$;\\
\item $\displaystyle\app3{a_1,a_2;\;b_1,b_2}{c}{x,\frac{x}{x-1}}$ and\\
$\displaystyle x^{1-c}\,(1-x)^{a_2}\,\hpg32{1+a_1+a_2-c,1+b_1+a_2-c,1-b_2}
{1+a_1+a_2+b_1-c,1+a_2-b_2}{1-x}$;\\
\item $\displaystyle\app3{1-b_1,1-b_2;\,b_1,b_2}{c}{x,\frac{x}{2x-1}}$ and\\
$\displaystyle(1-x)^{c-1}\,(1-2x)^{b_2}\,
\hpg21{\frac12(b_2-b_1+c),\frac12(b_1+b_2+c-1)}{c}{4x(1-x)}$;\\
\item $\displaystyle\app3{1-b_1,1-b_2;\,b_1,b_2}{3/2}
{-\frac{(t-1)^2}{4t},\frac{s(t-1)^2}{(1-s)(st^2-1)}}$ and\\
$\displaystyle t^{b_1}(1-t)^{-1}
\,\hpg21{1-b_2,b_2}{b_1+\frac12}{\frac{st^2}{st^2-1}}$. 
\end{itemize}
In the last pair, $s$ is any constant. The function pairs represent all cases
when a univariate specialization of Appell's $F_3$ function satisfies a second order 
ordinary Fuchsian equation fully implied by \mbox{\rm (\ref{app3a})--(\ref{app3b})}.
\end{theorem}
\proof Apart from the first trivial case, we employ correspondence (\ref{app2app3}) in each case of Theorems \ref{th:hpg32} and \ref{th:f2cases}. We substitute
\begin{eqnarray*}
x\mapsto \frac1x,\quad y\mapsto \frac1y,\quad c_1\mapsto 1+b_1-a_1, 
\quad c_2\mapsto 1+b_2-a_2, \quad a\mapsto 1+b_1+b_2-c,
\end{eqnarray*}
get rid of the power factors $x^{b_1}y^{b_2}$, and choose sometimes a more convenient companion
from the 24 Kummer's or six $\hpgo32$ series of the second function.\qed\\

In the first and fourth function pairs, a direct identity between the two functions holds in a neighborhood 
of $x=0$. In the second pair, the following formula around $x=0$ holds if $\mbox{Re}(c-a_2-b_2)>0$:
\begin{equation}
\app3{a_1,a_2;\;b_1,b_2}{c}{x,\,1}=
\frac{\Gamma(c)\Gamma(c-a_2-b_2)}{\Gamma(c-a_2)\Gamma(c-b_2)}
\hpg32{a_1,b_1,c-a_2-b_2}{c-a_2,c-b_2}{x}.
\end{equation}
Different from (\ref{eq:wrongformula}), here convergence of the $\hpgo21(1)$ coefficients
to powers of $x$ in the $F_3(x,1)$ expansion actually improves with higher powers of $x$; see
\cite{Vidunas09}.
In the third and last cases, relations between Appell's $F_3$ and $\hpgo21$ solutions of the implied ordinary differential equations can be derived following their analytic continuation.

The $F_3(x,x/(x-1))$ and $F_3(x,x/(2x-1))$ functions are considered in \cite{Karlsson80}. 
Karlsson expresses $F_3(x,x/(x-1))$ functions with $c=a_1+a_2+b_1+b_2$ or $c=a_1+b_2$ in terms
of Gauss hypergeometric functions; these relations can be obtained by simplifying the respective
$\hpgo32$ function in Theorem \ref{th:f3cases}. Karlsson's relation between $F_3(x,x/(2x-1))$ and 
$\hpgo21(4x(1-x))$ functions differs from ours by Euler's transformation \cite[(2.2.7)]{AAR99}.

Translating the case (\ref{eq:f2last})--(\ref{eq:f2hpgs}) of variable separation, we conclude
that the bivariate functions
\begin{equation} \label{eq:f3last}
\app3{1-b_1,1-b_2;\;b_1,b_2}{3/2}{-\frac{(t-s)^2}{4st},-\frac{(t-s)^2}{(s^2-1)(t^2-1)}}
\end{equation}
and
\begin{equation} \label{eq:f3hpgs}
\frac{s^{b_1}t^{b_1}(1-s^2)^{b_2}(1-t^2)^{b_2}}{s-t}\hpg21{b_1+b_2-\frac12,b_2}{b_1+\frac12}{s^2}
\hpg21{b_1+b_2-\frac12,b_2}{b_1+\frac12}{t^2}
\end{equation}
satisfy the same system of partial differential equations with respect to $s$ and $t$.

\section{Identities with Appell's $F_4$ function}
\label{sec:appellf4}

A system partial differential equations for Appell's $F_4$ function is:
\begin{eqnarray} \label{app4a}
x(1-x)\frac{\partial^2F}{\partial x^2}-y^2\frac{\partial^2F}{\partial y^2}
-2xy\frac{\partial^2F}{\partial x\partial y}+\left(c_1-(a+b+1)x\right)\frac{\partial F}{\partial x}
\nonumber\\ -(a+b+1)y\frac{\partial F}{\partial y} -abF\equal 0,\\ 
\label{app4b} y(1-y)\frac{\partial^2F}{\partial y^2}-x^2\frac{\partial^2F}{\partial x^2}
-2xy\frac{\partial^2F}{\partial x\partial y}+\left(c_2-(a+b+1)y\right)\frac{\partial F}{\partial y}
\nonumber \\ -(a+b+1)x\frac{\partial F}{\partial x}-abF\equal 0.
\end{eqnarray}
The singular locus on $\PP^2$ is $x=0$, $y=0$, the line at infinity, and 
\begin{equation} \label{app4sing}
x^2+y^2+1=2xy+2x+2y.
\end{equation}
The quadratic singular locus can be loosely described by the equation $\sqrt{x}+\sqrt{y}=1$. 
It can be parametrized as follows:
\begin{equation} \label{app4para}
x=t^2,\qquad y=(1-t)^2.
\end{equation}

First we characterize univariate specializations outside the singularity curves. 
The general case can be identified as Bailey's case (\ref{eq:bailey}) of variable separation.
\begin{theorem}
Suppose that a univariate specialization of Appell's $F_4(x,y)$ function 
is not onto singular locus $xy=0$ or $(\ref{app4sing})$, and satisfies a second order 
ordinary differential equation fully implied by \mbox{\rm(\ref{app4a})--(\ref{app4b})}.
Then the specialization is  represented by the general expression 
\begin{eqnarray} \label{app4:family}
\app4{a;\;b}{c,a+b-c+1}{s\,t,\,(1-s)(1-t)},
\end{eqnarray}
where $s$ is a constant.
The general second order equation is satisfied by $\displaystyle\hpg21{a, b}{c}{t}$.
\end{theorem}
\proof  The rank of the differential system (\ref{app4a})--(\ref{app4b}) is 4, and 
the partial differential equations of order 2 are linearly generated by (\ref{app4a})--(\ref{app4b}).
Linear elimination of $\partial^2 F/\partial x\partial y$, $\partial^2 F/\partial y^2$ 
from expression (\ref{eq:f2d2g}) for $d^2F/dt^2$ 
gives the following coefficient to $\partial^2 F/\partial x^2$:
\begin{equation} \label{eq:genf4}
\dot{x}^{2}+\frac{1-x-y}{y}\,\dot{x}\dot{y} +\frac{x}{y}\,\dot{y}^{2}.
\end{equation}
This expression must be equated to $0$. 
The expression is a quadratic form in the derivatives, 
with the discriminant equal to $4(x^2-2xy+y^2-2x-2y+1)/y^2$. 
We can factorize (\ref{eq:genf4}) into linear differential forms 
if we parametrize the surface $w^2=x^2-2xy+y^2-2x-2y+1$. 
A parametrization is easy if we set $w=y+u$; then
\begin{equation} \label{eq:solvy4}
y=\frac{(x+u-1)(x-u-1)}{2(x+u+1)}.
\end{equation}
Expression (\ref{eq:genf4}) factorizes as follows:
\begin{equation}  \label{eq:facteq4}
\frac{(x^2+2xu+u^2-2x+2u+1)^2}{2(x+u+1)^3(x+u-1)(x-u-1)}
(du+dx)\left(xdu-(u+1)dx\right),
\end{equation}
where we write $du$, $dx$ instead of $du/dt$, $dx/dt$ because the variable $t$ is irrelevant. 
The non-differential numerator factor defines the singular locus, since its parametrization
$x=t^2$, $u=-(t-1)^2$ is translated to (\ref{app4para}) by (\ref{eq:solvy4}).
The differential factors give the following two families of solutions:
\begin{equation}
u=C-x, \qquad  u=Cx-1.
\end{equation}
They translate into the following relations between $x$ and $y$:
\begin{equation}
y=\frac{C-1}{C+1}x+\frac{1-C}2, \qquad  y=-\frac{1-C}2x+\frac{C-1}{C+1}.
\end{equation}
The two solution families actually coincide, since they are related by
the parameter transformation $C\mapsto(3-C)/(1+C)$. After the substitutions
$C\mapsto 2s-1$, $x\mapsto st$, the first family is described by $y=(1-s)(1-t)$, 
giving us (\ref{app4:family}).

Adding equation (\ref{app4a}) multiplied by $s$ with equation (\ref{app4b}) multiplied by $1-s$,
and the substitutions $x=st$, $y=(s-1)(t-1)$ gives us
\begin{eqnarray*}
t(1-t)\left(s^2\frac{\partial^2 F}{\partial x^2}+2s(s-1)\frac{\partial^2 F}{\partial x\partial y}
+(s-1)^2\frac{\partial^2 F}{\partial y^2}\right)+\left(c_1-(a+b+1)t\right)s\frac{\partial F}{\partial x}\\
+\left(a+b+1-c_2-(a+b+1)t\right)(s-1)\frac{\partial F}{\partial y}-abF \equal 0.
\end{eqnarray*}
If $c_2=a+b+1-c_1$, we recognize Euler's equation (\ref{eq:euler}) in $t$.
\qed\\

For fixed $s$, the arguments $x$ and $y$ in (\ref{app4:family}) are linearly related. 
The linear relations are precisely the tangent lines to the singularity curve (\ref{app4sing}). 
If we let both $s$ and $t$ vary, their symmetry 
indicates Bailey's case \cite{Bailey33} of variable separation  
for differential system (\ref{app4a})--(\ref{app4b}) when $c_1+c_2=a+b+1$.
The conclusion is that 
\begin{eqnarray} \label{app4:separate}
\app4{a;\;b}{c,a+b-c+1}{st,(1-s)(1-t)}\quad\mbox{and}\quad
\hpg21{a, b}{c}{s}\hpg21{a, b}{c}{t}
\end{eqnarray}
satisfy the same system of partial differential equations (with respect to $s$ and $t$).
This is illustrated by identity (\ref{eq:bailey}). 

The $F_2(x,y)$ case of variable separation in (\ref{eq:f2last})--(\ref{eq:f2hpgs}) 
is related to Bailey's case via the following identity, with $c=2b$:
\begin{equation} \label{eq:f2f4rel}
\app4{a;\;b}{c,a-b+1}{x,y^2}=(1+y)^{-2a}\app2{a; b, a-b+\frac12}{c,2a-2b+1}
{\frac{x}{(1+y)^2},\frac{4y}{(1+y)^2}}.
\end{equation}
This identity is derived by Srivastava in \cite{Srivastava66}; 
it is presented in \cite[9.4.(215)]{Srivastava85}. An equivalent identity is obtained
in \cite[pg. 27]{Appel26}, and presented in \cite[9.4.(97)]{Srivastava85}.
The Appell functions in (\ref{eq:f2last}) and (\ref{app4:family}) are generic $F_2$ and $F_4$ functions
satisfying the quadratic condition of Sasaki--Yoshida \cite[Sections 4.3, 5.5]{SasakiYoshida}:
four linearly independent solutions of their partial differential equations systems are quadratically related. 

In \cite[Section 5.4]{SasakiYoshida}, the following identity between $F_4$ and $F_2$
functions is implied: 
\begin{equation} \label{eq:f2f4rel2}
\app4{\frac{a}2;\frac{a+1}2}{b_1\!+\!\frac12,b_2\!+\!\frac12}{x^2,y^2}=
(1+x+y)^{-a}\app2{a;b_1,b_2}{2b_1,2b_2}{\frac{2x}{x+y+1},\frac{2y}{x+y+1}}.
\end{equation}
Equivalent identities are presented in \cite{Bailey38}, \cite[9.4.(175)]{Srivastava85} and 
\cite[Lemma 5.2]{Kato00}. This identity has no application for $F_2(x,2-x)$ functions. 

There remains to consider specialization of the $F_4$ function to the quadratic singularity
locus (\ref{app4sing}). Due to parametrization (\ref{app4para}), we are looking at
$F_4(t^2,(t-1)^2)$ 
functions. Can such a function satisfy a second order Fuchsian equation? 
\begin{theorem} \label{th:app4ga}
The univariate functions
\begin{equation}
\app4{a;\;b}{c,a+b-c+\frac32}{t^2,(1-t)^2} \qquad\mbox{and}\qquad
\hpg21{2a,\,2b}{2c-1}{t}
\end{equation}
satisfy the same Fuchsian equation. This is the only case when the differential system 
$(\ref{app4a})$--$(\ref{app4b})$ fully implies a second order Fuchsian equation.
\end{theorem}
\proof Let $J$ be the differential ideal generated by (\ref{app4a})--(\ref{app4b}). 
After substitution (\ref{app4para}) into the coefficients, the rank of partial derivatives becomes 3,
and there is another independent linear relation between the partial derivatives of order at most 2.
Yet, second order partial derivatives are eliminated from the full differential 
\begin{eqnarray} \label{ap4:sdervs2}
\frac{d^2F}{dt^2}\equal4t^2\frac{\partial^2F}{\partial x^2}+8t(t-1)\frac{\partial^2 F}{\partial x\partial y}
+4(t-1)^2\frac{\partial^2F}{\partial y}+2\frac{\partial F}{\partial x}+2\frac{\partial F}{\partial y}
\end{eqnarray}
already by equations (\ref{app4a})--(\ref{app4b}).
After we substitute (\ref{app4para}) into the coefficients in \mbox{(\ref{app4a})--(\ref{app4b})}, 
add the first equation multiplied by $4/(1-t)$ and the second equation multiplied by $4/t$, we get
\begin{eqnarray}  \label{ap4:ddxy}
4t^2\frac{\partial^2F}{\partial x^2}+8t(t-1)\frac{\partial^2F}{\partial x\partial y}
+4(t-1)^2\frac{\partial^2F}{\partial y^2}+4\frac{(a+b+1)t-c_1}{t-1}\frac{\partial F}{\partial x}\nonumber\\ 
+4\frac{(a\!+\!b\!+\!1)t+c_2\!-\!a\!-\!b\!-\!1}{t}\frac{\partial F}{\partial y}+\frac{4ab}{t(t-1)}F\equal0.
\end{eqnarray}
Identifying (\ref{ap4:sdervs2}) and
\begin{eqnarray} \label{ap4:sdervs}
\frac{dF}{dt}\equal2t\frac{\partial F}{\partial x}+2(t-1)\frac{\partial F}{\partial y}
\end{eqnarray}
we recognize in (\ref{ap4:ddxy}) the expression
\begin{eqnarray} \label{ap4:ddty}
\frac{d^2F}{dt^2}+\frac{(2a\!+\!2b\!+\!1)t+1\!-\!2c_1}{t(t-1)}\frac{dF}{dt}
+\frac{4c_1\!+\!4c_2\!-\!4a\!-\!4b\!-\!6}{t}\frac{\partial F}{\partial y}+\frac{4ab}{t(t-1)}F=0.
\end{eqnarray}
We have differentiation here with respect to $t$ only if $c_1+c_2=a+b+\frac32$. 
Euler's equation is recognizable.
\qed\\

The general ordinary differential equation for $\displaystyle\app4{a;\;b}{c_1,c_2}{t^2,(t-1)^2}$
is computed  by Kato \cite[(4.1)]{Kato95}. We reproduce the differential equation in 
Theorem \ref{th:a4qsing} here below as demonstration of general computational techniques;
see formula (\ref{ap4:singeq}) in Theorem \ref{th:a4qsing}.
The most important knowledge at this point is that the differential equation is Fuchsian
with singularities at $t=0$, $t=1$, $t=\infty$, and the following local exponents:
\begin{eqnarray} \label{eq:apqle}
\begin{array}{llcl}
\mbox{at $t=0$}: & 0, & 2-2c_1, & c_2-a-b;\\
\mbox{at $t=1$}: & 0, & 2-2c_2, & c_1-a-b;\\
\mbox{at $t=\infty$}: & 2a, & 2b, & c_1+c_2-1.
\end{array}
\end{eqnarray}
If some difference of local exponents at the same singular point is equal to 1, 
then there is a chance that the accessory parameter is right and the Fuchsian equation
is projectively equivalent to Appendix equation (\ref{eq:hpgde3}) for the $\hpgo32$ function.
We present these cases in the following theorem.

Incidentally, we recall some Kato's results on reducibility of the monodromy groups of differential system 
(\ref{app4a})--(\ref{app4b}) or the ordinary equation for $F_4(t^2,(1-t)^2)$. The monodromy group of differential system (\ref{app4a})--(\ref{app4b}) with $c_1,c_2\not\in\ZZ$ is reducible if and only if 
one of the numbers
\begin{equation} \label{eq:ap4red}
a,\; b,\;  c_1-a,\; c_1-b,\; c_2-a,\; c_2-b,\; c_1+c_2-a,\; c_1+c_2-b
\end{equation}
is an integer  \cite[Section 8]{Kato95}. If no difference of the local exponents in (\ref{eq:apqle}) at a singular point is an integer, then Kato's ordinary equation is reducible if and only if 
one of the numbers in (\ref{eq:ap4red}) or $c_1+c_2-a-b-\frac12$
is an integer \cite[Section 14]{Kato95}.
\begin{theorem}
The univariate functions 
\begin{equation} \label{eq:a4hpg3}
\app4{a;\,b}{c+\frac12,\frac12}{t^2,(1-t)^2} \qquad\mbox{and}\qquad 
\hpg32{2a,\,2b,\,c}{a+b+\frac12,2c}{t}
\end{equation}
satisfy the same ordinary Fuchsian equation.
That equation is satisfied by
\begin{eqnarray} \label{eq:a4hpg3a}
&&(1-t)\;\app4{a+\frac12;\,b+\frac12}{c+\frac12,\frac32}{t^2,(1-t)^2},\\ \label{eq:a4hpg3b}
&&(1-t)^{-2a}\,\app4{a;\,a+\frac12}{c+\frac12,1+a-b}{\frac{t^2}{(t-1)^2},\frac1{(t-1)^2}}
\end{eqnarray}
as well. This represents all cases with $c_1,c_2\not\in\ZZ$ when the monodromy of 
\mbox{\rm (\ref{app4a})--(\ref{app4b})} is irreducible and Kato's differential equation 
$(\ref{ap4:singeq})$ is projectively equivalent to hypergeometric equation $(\ref{eq:hpgde3})$.
\end{theorem}
\proof Here are all cases when a difference of local exponents is equal to $1$:
\begin{itemize}
\item $c_2=\frac12$. Once we identify the local exponent sets (\ref{eq:apqle}) and (\ref{eq:hpg3le})
under this condition, Kato's equation coincides with (\ref{eq:hpgde3}) up to variable notation and a polynomial factor. We get the relation between the functions in (\ref{eq:a4hpg3}) after the 
substitution \mbox{$c_1\mapsto c_1+\frac12$} (for a better form). 
Functions (\ref{eq:a4hpg3a})--(\ref{eq:a4hpg3b}) are solutions of the same ordinary equation by Appendix transformations (\ref{eq:a4sola})--(\ref{eq:a4solz}).
\item $c_2=\frac32$. To identify the local exponents in (\ref{eq:apqle}) and (\ref{eq:hpg3le}), 
we apply the projective transformation $y(t)\mapsto (1-t)^{-1}y(t)$. This shifts the local exponents
at $t=1$ by $1$, and the local exponents at $t=\infty$ by $-1$. We get the coincidence of Kato's and hypergeometric equations unconditionally. A direct conclusion is that the functions
\begin{eqnarray*} 
\app4{a;\,b\,}{c,\,\frac32}{t^2,(1-t)^2} \qquad\mbox{and}\qquad 
\frac1{1-t}\;\hpg32{2a-1,\,2b-1,\,c-\frac12}{a+b-\frac12,2c-1}{t}
\end{eqnarray*}
satisfy the same ordinary Fuchsian equation. After multiplication of both function by $1-t$ and a
simple transformation of parameters we get the $\hpgo32$ function as in (\ref{eq:a4hpg3}),
and the $F_4$ function as in (\ref{eq:a4hpg3a}).
\item $c_1=a+b+1$. The local exponents in (\ref{eq:apqle}) and (\ref{eq:hpg3le}) can be identified directly. The accessory parameter is right if $c_2=\frac12$ or $a=0$ or $b=0$. In the later two cases,
the monodromy group of the system (\ref{app4a})--(\ref{app4b}) is reducible.
\item $c_1=a+b-1$. The local exponents in (\ref{eq:apqle}) and (\ref{eq:hpg3le}) can be identified 
after the projective transformation $y(t)\mapsto (1-t)^{-1}y(t)$. The accessory parameter is right
if $c_2=\frac32$ or $a=1$ or $b=1$. The later two cases are reducible.
\item $c_1+2c_2=a+b+1$. The local exponents in (\ref{eq:apqle}) and (\ref{eq:hpg3le}) can be identified 
after the projective transformation $y(t)\mapsto (1-t)^{1-2c_2}y(t)$. The accessory parameter is right
if $c_2=\frac12$ or $c_2=a$ or $c_2=b$. The later two cases are reducible.
\item $c_1+2c_2=a+b+3$. The local exponents in (\ref{eq:apqle}) and (\ref{eq:hpg3le}) can be identified 
after the projective transformation $y(t)\mapsto (1-t)^{2-2c_2}y(t)$. The accessory parameter is right
if $c_2=\frac32$ or $c_2=a+1$ or $c_2=b+1$. The later two cases are reducible.
\item $b=a+\frac12$. To identify the local exponents in (\ref{eq:apqle}) and (\ref{eq:hpg3le}), 
we permute the singular points $t=1$ and $t=\infty$, and shift the local exponents there by $\pm2a$.
Then Kato's and hypergeometric equations coincide unconditionally. 
A direct conclusion is that the functions
\begin{eqnarray*} 
\app4{a,a+\frac12}{c_1,c_2}{t^2,(1-t)^2} \quad\mbox{and}\quad
(1-t)^{2a}\,\hpg32{2a,\,2a-2c_2+1,c_1-\frac12}{2a-c_2+1,2c_1-1}{\frac{t}{t-1}}
\end{eqnarray*}
satisfy the same ordinary Fuchsian equation. After straightforward transformations 
we get the $\hpgo32$ function as in (\ref{eq:a4hpg3}),
and the $F_4$ function as in (\ref{eq:a4hpg3b}).
\item $c_1+c_2=2a$. To identify the local exponents in (\ref{eq:apqle}) and (\ref{eq:hpg3le}),
we permute the singular points $t=1$ and $t=\infty$, and shift the local exponents there by $\pm(2a-1)$.
The accessory parameter is right if $a=b+\frac12$ or $a=1$ or $c_1=a$. The later two cases are reducible.
\item $c_1+c_2=2a+2$. To identify the local exponents in (\ref{eq:apqle}) and (\ref{eq:hpg3le}),
we permute the singular points $t=1$ and $t=\infty$, and shift the local exponents there by $\pm2a$.
The accessory parameter is right if $b=a+\frac12$ or $a=0$ or $c_1=a+1$. 
The later two cases are reducible.
\item The cases $c_1=\frac12$, $c_1=\frac32$, $c_2=a+b+1$, $c_2=a+b-1$, 
$2c_1+c_2=a+b+1$, $2c_1+c_2=a+b+3$, $a=b+\frac12$, $c_1+c_2=2b$, $c_1+c_2=2b+2$
are symmetric to the above.
\end{itemize}
Non-reducible cases fall into the cases when $c_1$ or $c_2$ is equal to $\frac12$ or $\frac32$, 
or $a-b=\pm\frac12$.
\qed\\

Similarly, we may find the cases when Kato's differential equation is a symmetric tensor square of Euler's hypergeometric equation (\ref{eq:euler}). Recall that if $y_1$, $y_2$ is a basis of second order linear ordinary differential equation, then $y_1^2$, $y_2^2$, $y_1y_2$ form a basis of the symmetric square 
differential equation. One case when Kato's differential equation is a symmetric square of (\ref{eq:euler}) follows from Bailey's case of variable separation. If we set $s=t$ in (\ref{app4:separate}), we conclude that the functions
\begin{equation} \label{eq:f4sep2}
\app4{a;\;b}{c,a+b-c+1}{t^2,(1-t)^2}
\qquad \mbox{and} \qquad \hpg21{a,b\,}{c}{\,t}^2
\end{equation}
satisfy the same ordinary Fuchsian equation. Appell himself derived this fact explicitly \cite{Appell84}.
Here is characterization of other such cases.

\begin{theorem}
The univariate functions 
\begin{equation} \label{eq:a4hpg2s}
\app4{2c-\frac12;\,3c-1}{c+\frac12,c+\frac12}{t^2,(1-t)^2} \qquad\mbox{and}\qquad 
\hpg21{c,\,3c-1}{2c}{t}^2
\end{equation}
satisfy the same ordinary Fuchsian equation of order $3$.
That equation is satisfied by
\begin{eqnarray} \label{eq:a4hpg2sa}
&&t^{1-2c}\,\app4{c;\,2c-\frac12}{\frac32-c,c+\frac12}{t^2,(1-t)^2},\\ \label{eq:a4hpg2sb}
&&t^{1-2c}(1-t)^{1-2c}\,\app4{\frac12;\,c}{\frac32-c,\frac32-c}{t^2,(1-t)^2}.
\end{eqnarray}
as well. Together with $(\ref{eq:f4sep2})$ this represents all 
cases when Kato's differential equation $(\ref{ap4:singeq})$ is projectively equivalent to 
a symmetric tensor square of Euler's hypergeometric equation $(\ref{eq:euler})$.
\end{theorem}
\proof A necessary condition on local exponents (\ref{eq:apqle}) is that at each singular point
one  local exponent is the arithmetic mean of the other two. If we pick for the arithmetic mean
at $t=0$ the exponent $c_2-a-b$, or at $t=1$ the exponent $c_1-a-b$, 
or at $t=\infty$ the exponent $c_1+c_2-1$, we arrive at Appell's case (\ref{eq:f4sep2}).
Up to symmetries of parameters, there are three different ways to pick up other local exponents
as arithmetic means:
\begin{itemize}
\item The exponents $2-2c_1$, $2-2c_2$, $2a$ are the arithmetic means at 
$t=0$, $t=1$, $t=\infty$, respectively. We have 3 linear equations for the 4 parameters, 
and the general solution is presented in (\ref{eq:a4hpg2s}). The local exponents of the $\hpgo21$ function must be equal to the non-arithmetic mean exponents divided by 2. Explicit computation
confirms coincidence of ordinary differential equations; in other words, the accessory parameter turns out to be right. Functions (\ref{eq:a4hpg2sa})--(\ref{eq:a4hpg2sb}) are solutions of the same ordinary equation by Appendix transformations (\ref{eq:a4sola})--(\ref{eq:a4solz}).
\item The exponents $2-2c_1$, $0$, $2a$ are the arithmetic means.  A general parametrization of the parameters is the same as in (\ref{eq:a4hpg2sa}). Evidently, we get a related case.
\item The exponents $0$, $0$, $2a$ are the arithmetic means. A general parametrization of the parameters is the same as in (\ref{eq:a4hpg2sb}).  
\qed\\
\end{itemize}


\section{Identities with Appell's $F_1$ function}

A system partial differential equations for Appell's $F_1$ function is:
\begin{eqnarray} \label{app1a}
x(1-x)\frac{\partial^2F}{\partial x^2}+y(1-x)\frac{\partial^2F}{\partial x\partial y}
+\left(c-(a+b_1+1)x\right)\frac{\partial F}{\partial x}-b_1y\frac{\partial F}{\partial y}
-ab_1F=0,\\ \label{app1b}
y(1-y)\frac{\partial^2F}{\partial y^2}+x(1-y)\frac{\partial^2F}{\partial x\partial y}
+\left(c-(a+b_2+1)y\right)\frac{\partial F}{\partial y}-b_2x\frac{\partial F}{\partial x}
-ab_2F=0.
\end{eqnarray}
This is an equation system of rank 3 generally. 
The equation
\begin{equation}
\left((1-y)\frac{\partial F}{\partial y}-b_2\right)\mbox{[Eq. (\ref{app1a})]}
+\left((x-1)\frac{\partial F}{\partial x}+b_1\right)\mbox{[Eq. (\ref{app1b})]}
\end{equation}
can be computed to be $a-c+1$ times
\begin{equation} \label{eq:app1c}
(y-x)\frac{\partial^2F}{\partial x\partial y}
+b_2\frac{\partial F}{\partial x}-b_1\frac{\partial F}{\partial y}=0
\end{equation}
This equation is of order 2 as well; it holds even for the $F_1$ functions with $c=a+1$. 
In the case $c=a+1$, the equation system (\ref{app1a})--(\ref{app1b}) is not holonomic;
equation (\ref{eq:app1c}) has to be added to get a finite dimensional space of solutions.

The singular locus of this equation on $\PP^1\times\PP^1$ is the union of the following lines:
\begin{equation} \label{eq:singf1}
x=0,\quad x=1,\quad x=\infty, \quad y=0,\quad y=1,\quad y=\infty, \quad y=x.
\end{equation}
Univariate specializations to the singular lines are known since \cite{Appel26}:
\begin{eqnarray}
\app1{a;\;b_1,b_2}{c}{x,\,0}\equal\hpg21{a,\,b_1}{c}{x},\\
\app1{a;\;b_1,b_2}{c}{x,\,1}\equal\frac{\Gamma(c)\Gamma(c-a-b_2)}{\Gamma(c-a)\Gamma(c-b_2)}
\hpg21{a,\,b_1}{c-b_2}{x},\\
\app1{a;\;b_1,b_2}{c}{x,\,x}\equal\hpg21{a,\,b_1+b_2}{c}{x}.
\end{eqnarray}
The second identity holds if Re$(c-a-b_2)>0$.

\begin{theorem} \label{th:f1usde}
Suppose that a univariate specialization $F_1(x(t),y(t))$ is not onto any of singularity curves in 
$(\ref{eq:singf1})$. Then equations $(\ref{app1a})$--$(\ref{app1b})$ and $(\ref{eq:app1c})$ fully imply
a second order differential equation for $F_1(x(t),y(t))$ if and only if the functions $x(t)$ and $y(t)$ are related by the non-linear differential equation
\begin{eqnarray} \label{eq:f1xy2rel}
\frac{\ddot{x}}{\dot{x}}-\frac{\ddot{y}}{\dot{y}}
-(a+1)\left(\frac{\dot{x}}{x-1}-\frac{\dot{y}}{y-1}\right)
+c\left(\frac{\dot{x}}{x(x-1)}-\frac{\dot{y}}{y(y-1)}\right) \nonumber\\
+\frac{x(x-1)y(y-1)}{(y-x)\dot{x}\dot{y}}
\left(\frac{\dot{x}}{x}-\frac{\dot{y}}{y}\right)\!
\left(\frac{\dot{x}}{x(x-1)}-\frac{\dot{y}}{y(y-1)}\right)\!
\left(\frac{b_1\dot{x}}{x-1}+\frac{b_2\dot{y}}{y-1}\right)\equal 0.
\end{eqnarray}
\end{theorem}
\proof The second order derivative $d^2F/dt^2$ in (\ref{eq:f2d2g}) can be reduced by equations
\mbox{(\ref{app1a})--(\ref{app1b})} and (\ref{eq:app1c}) to first order partial derivatives without additional restrictions on $x(t),y(t)$. From equation
\begin{eqnarray*}
\frac{\dot{x}^2}{x(1-x)}\mbox{[Eq. (\ref{app1a})]}+\frac{\dot{y}^2}{y(1-y)}\mbox{[Eq. (\ref{app1b})]}
-\frac{(\dot{x}y-x\dot{y})^2}{xy\,(y-x)}\mbox{[Eq. (\ref{eq:app1c})]}
\end{eqnarray*}
we straightforwardly derive
\begin{eqnarray*}
\frac{d^2F}{dt^2}\equal
\left(\ddot{x}-\frac{c-(a+b_1+1)x}{x(1-x)}\dot{x}^2
+\frac{b_2\dot{x}^2}{x}+\frac{b_2(\dot{x}-\dot{y})^2}{y-x}+\frac{b_2(x+y-1)\dot{y}^2}{y(1-y)}\right)
\frac{\partial F}{\partial x}\\
&&+\left(\ddot{y}-\frac{c-(a+b_2+1)y}{y(1-y)}\dot{y}^2
+\frac{b_1\dot{y}^2}{y}-\frac{b_1(\dot{x}-\dot{y})^2}{y-x}+\frac{b_1(x+y-1)\dot{x}^2}{x(1-x)}\right)
\frac{\partial F}{\partial y}\\
&&+\frac{ab_1\dot{x}^2}{x(1-x)}+\frac{ab_2\dot{y}^2}{y(1-y)}.
\end{eqnarray*}
We can recognize the derivative $dF/dt$ on the right-hand side precisely when
the quotient of the coefficients to  $\partial F/\partial x$ and $\partial F/\partial y$ is
equal to $\dot{x}/\dot{y}$, according to  (\ref{eq:f2d1g}). This leads to equation (\ref{eq:f1xy2rel}).
\qed\\

Without loss of generality, we can set $\dot{x}=1$, $\ddot{x}=0$ in equation (\ref{eq:f1xy2rel}).
Then the non-linear equation can be written in the form $\ddot{y}=g_3(\dot{y})$, where $g_3(\dot{y})$
is a cubic polynomial in $\dot{y}$ with coefficients rational functions in $x,y$. This is not a canonical form of a non-linear differential equation 
without movable essential singularities \cite[Chapter XIV]{Ince56}.
Nevertheless, its Liouville invariants \cite[(14)--(15)]{Hietarinta02} are zero,
hence the equation is a differentiable transformation $x=\Phi(X,Y)$, $y=\Psi(X,Y)$ of the simplest equation $Y''=0$. However, it is not straightforward to find the transformation.

Using the identity
\begin{equation} \label{eq:f1f3}
\app1{a; b_1, b_2}{c}{x,y} = (1-x)^{-b_1}\app3{c-a,a;b_1,b_2}{c}{\frac{x}{x-1},y}.
\end{equation}
one can translate the two non-singular cases of Theorem \ref{th:f3cases} to a few univariate specializations of the $F_1$ function relevant to Theorem \ref{th:f1usde}. 
The $F_2(x,2-x)$ or $F_3(x,x/(2x-1))$ case can be translated as follows.
\begin{theorem}
The following two univariate functions satisfy the same ordinary Fuchsian differential equation
of order $2$:
\begin{equation} \label{eq:f1q0}
\app1{a;\;2b, a-b}{1+b}{x,\,x^2} \qquad\mbox{and}\qquad
(1-x)^{-2a}\,\hpg21{a,\;\frac12}{1+b}{-\frac{4x}{(x-1)^2}}.
\end{equation}
The same equation is satisfied by
\begin{eqnarray} \label{eq:f1q1}
(1\!-\!x)^{-a}\app1{\!a; 1\!-\!a, a\!-\!b}{1+b}{\frac{x}{x-1}, -x},\
\left(1\!-\!x^2\right)^{-a}\app1{\!a; 2b, 1\!-\!a}{1+b}{\frac{x}{x+1}, \frac{x^2}{x^2-1}},\\ \label{eq:f1q2}
\app1{a;\,2b, a-b}{2a}{1-x,\,1-x^2},\qquad\,
x^{-a}\app1{a;\, a-b, a-b}{2a}{\frac{x-1}x, 1-x},\\  \label{eq:f1q3}
x^{-a}\app1{a;\, a-b, a-b}{1+a-2b}{\frac1x,\,x},\quad
(x-1)^{-a}\app1{a;\, 1-a, a-b}{1+a-2b}{\frac1{1-x}, x+1},\\ \label{eq:f1q4}
x^{-b}(1-x)^{2b-2a}\app1{1-2b;\, a-b, a-b}{1+a-2b}{\frac1{1-x}, \frac{x}{x-1}},\\ \label{eq:f1q5}
x^{-b}(1-x)^{1-2a}\app1{1-2b;\, a-b, 1-a}{1+a-2b}{x+1,\,x}.
\end{eqnarray}
\end{theorem}
\proof The $F_3$ function in (\ref{eq:f1f3}) can be directly related to the penultimate case of Theorem 
\ref{th:f3cases} if $a=1-b_2$ and $c=2-b_1-b_2$. The direct conclusion is that the functions
\begin{equation} \label{eq:f1q6}
\app1{1-b_2;\;b_1,b_2}{2-b_1-b_2}{x,\,\frac{x}{x+1}} \quad\mbox{and}\quad 
\frac{(1+x)^{b_2}}{1-x}\,\hpg21{1-b_1,\;\frac12}{2-b_1-b_2}{-\frac{4x}{(x-1)^2}}
\end{equation}
satisfy the same second order ordinary Fuchsian equation. We rename $b_1\mapsto 1-a$, 
$b_2\mapsto a-b$, and adjust the power factor to the $\hpgo21$ function, so to get the second
function in (\ref{eq:f1q0}). As presented in Appendix formulas (\ref{eq:f1s1})--(\ref{eq:f1s9}), there are generally 60 $F_1$ bivariate series giving solutions of the differential system (\ref{app1a})--(\ref{app1b}).
Investigation of those 60 $F_1$ series under the given specialization 
gives the $F_1$ functions in (\ref{eq:f1q0})--(\ref{eq:f1q5}).
\qed\\

\noindent
The $F_1$ functions in (\ref{eq:f1q0})--(\ref{eq:f1q5}) represent all possible relations 
between the 2 variables and the parameters among the specialized 60 $F_1$ series, 
obtainable from Appendix formulas (\ref{eq:f1s1})--(\ref{eq:f1s9}) and formally giving 
univariate solutions of the same ordinary differential equation. In particular, the $F_1$ function
in (\ref{eq:f1q6}) has the same relation between the arguments and parameters (up to permutation
of the two arguments and the $b_1,b_2$ parameters) as the first function in (\ref{eq:f1q1}).
The four functions in (\ref{eq:f1q0})--(\ref{eq:f1q1}) can be identified in a neighborhood of $x=0$,
while the two functions in (\ref{eq:f1q2}) can be identified in a neighborhood of $x=1$. 
Notice that relation between the two arguments in (\ref{eq:f1q4}) or (\ref{eq:f1q5}) is linear.

The variable separation case (\ref{eq:f3last})--(\ref{eq:f3hpgs}) for the $F_3$ function translates
into a reducible transformed system (\ref{app1a})--(\ref{app1b}). 
A representative conclusion is that the differential system for
\begin{equation}
\app1{b+\frac12;\,b,\,\frac12-b}{3/2}{\frac{(s-t)^2}{(s+t)^2},-\frac{(s-t)^2}{(s^2-1)(t^2-1)}}
\end{equation}
has the following elementary solution:
\begin{equation}
\frac{1}{s-t}\,(s+t)^{2b}\left(1-s^2\right)^{\frac12-b}\left(1-t^2\right)^{\frac12-b}.
\end{equation}

If differential condition (\ref{eq:f1xy2rel}) of Theorem \ref{th:f1usde} is tried for linear solutions $y(x)$, 
the following cases are found:
\begin{itemize}
\item Reducible Appell's functions $\displaystyle\app1{a;b_1,b_2}{a+1}{x,sx}$,
$\displaystyle\app1{a;b_1,c-b_1}{c}{1-x,1-sx}$ 
and $\displaystyle\app1{a;0,b_2}{c}{x,s}$,  where $s$ is a constant.
\item The cases (\ref{eq:f1q4}), (\ref{eq:f1q5}) with $y=1-x$, $y=x-1$, respectively,
and a symmetric case with $y=x+1$.
\item The trivial identity $\displaystyle\app1{-1;\,b_1,b_2}{c}{x,\frac{c-b_1x}{b_2}}=0$.
\end{itemize}

\section{Computational aspects}
\label{sec:compute}

Identities between bivariate and univariate hypergeometric series are usually derived by methods of series manipulation, or using integral representations. The method of relating differential equations is usually considered as tedious and computationally costly. In \cite[pg. 314]{Srivastava85}
Srivastava and Karlsson characterize this method as follows:
\begin{quote}
This method (though simple in theory) is rather laborious in practice,
and is not very useful for discovering new transformations.
\end{quote}
They immediately mention that Appell himself \cite{Appell84} showed that the functions in
(\ref{eq:f4sep2}) satisfy the same differential equation of order 3.  

Here we demonstrate the basic computational routine to translate linear differential equations
in partial derivatives to ordinary differential equations for univariate specialization of their solutions.
In our terminology, we seek to recognize {\it partial differential forms} of ordinary equations within
{\em specialized forms} of given systems of partial differential equations. 

In the specialization setting $x\mapsto x(t)$, $y\mapsto y(t)$ of Definition \ref{def:whole}, let $J$ denote the left ideal of the Weyl algebra \mbox{$\CC[x,y]\langle\partial/\partial x,\partial/\partial y\rangle$} determined by a given system of partial differential equations. We work in the $\CC(t)\langle d/dt\rangle$ left module $M$ that is linearly generated over $\CC(t)$ by the full derivatives $1, d/dt, d^2/dt^2$, etc., and the partial derivatives $\partial/\partial x, \partial/\partial y, \partial^2/\partial x^2, \partial^2/\partial x\partial y$, etc. 
We let $d/dt$ act from the left following the Leibniz rule, with the action on the partial derivatives
defined by the identity 
\begin{equation} \label{eq:dtdxdy}
\frac{d}{dt} = \dot{x}\,\frac{\partial}{\partial x}+\dot{y}\,\frac{\partial}{\partial y},
\end{equation}
compatible with (\ref{eq:f2d1g}). Let $P$ be the submodule of $M$ generated by 
\mbox{$d/dt-\dot{x}\,\partial/\partial x-\dot{y}\,\partial/\partial y$}. The elements of $P$ give expressions
of the full derivatives $d/dt,d^2/dt^2$, etc., in terms of partial derivatives, like the series of formulas starting with (\ref{eq:f2d1g})--(\ref{eq:f2d2g}).

Let $\widetilde{J}$ denote the submodule of $M$ generated by special forms of the elements in $J$.
Let $\widetilde{P}$ denote the union of $\widetilde{J}$ and $P$, as a submodule of $M$. 
We propose to use Gr\"obner basis computations for $\widetilde{P}$ with respect to an ordering of
full and partial derivatives suitable for elimination of the partial derivatives. 

Here we briefly demonstrate
this strategy by computing Kato's differential equation \cite[(4.1)]{Kato95} for the general 
$F_4(t^2,(1-t)^2)$ function. We used this equation in Section \ref{sec:appellf4}, and presented
its singularities and local exponents in (\ref{eq:apqle}).
\begin{theorem} \label{th:a4qsing}
The univariate function $\displaystyle\app4{a;\;b}{c_1,c_2}{t^2,(t-1)^2}$ satisfies the differential equation
\begin{eqnarray}  \label{ap4:singeq}
\frac{d^3F}{dt^3}
+\left(\frac{a+b+2c_1-c_2+1}t+\frac{a+b-c_1+2c_2+1}{t-1}\right)\frac{d^2F}{dt^2}
\qquad\qquad\nonumber\\
+\left(\frac{(2c_1-1)(a+b-c_2+1)}{t^2}+2\frac{a+b-c_1-c_2+1+2ab+2c_1c_2}{t(t-1)}\right.
\qquad\nonumber\\
+\left.\frac{(2c_2-1)(a+b-c_1+1)}{(t-1)^2}\right)\!\frac{dF}{dt}
+\frac{2ab\left(2(c_1\!+\!c_2\!-\!1)t-2c_1+1\right)}{t^2(t-1)^2}F \equal0.
\end{eqnarray} 
\end{theorem}
\proof In the course of the proof of Theorem \ref{th:app4ga} we basically derived equation (\ref{ap4:ddty}) 
as an element of $\widetilde{P}$. We differentiate (\ref{ap4:ddty}) with respect to $t$, keeping in mind that $d/dt$ acts on the partial derivative $\partial F/\partial y$ following (\ref{eq:dtdxdy}), or more specifically, following (\ref{ap4:sdervs}). We get the partial derivatives $\partial^2F/\partial x\partial y$
and $\partial^2F/\partial y^2$ in the computation; we can eliminate them using original equations
(\ref{app4a})--(\ref{app4b}), but then first order derivatives $\partial F/\partial x$, $\partial F/\partial y$ occur. For linear elimination of these derivatives, 
we use (\ref{ap4:sdervs}) and the equation (\ref{ap4:ddty})  itself.
Equation (\ref{ap4:singeq}) is obtained as
\begin{eqnarray*}
&&\hspace{-12pt} \frac{d}{dt}\mbox{[Eq. (\ref{ap4:ddty})]}
+\frac{2(2c_1+2c_2-2a-2b-3)}{(t-1)^2}\left(\mbox{[Eq. (\ref{app4a})]}
+\frac{1-t^2}{t^2}\mbox{[Eq. (\ref{app4b})]}\right)\\
&&\hspace{-12pt}
+\frac{(c_1\!-\!a\!-\!b\!-\!1)(2c_1\!+\!2c_2\!-\!2a\!-\!2b\!-\!3)}{t(t-1)^2}\mbox{[Eq.~(\ref{ap4:sdervs})]}
+\frac{(c_1\!+\!c_2\!+\!1)t \!+\! c_2\!-\!a\!-\!b\!-\!2}{t(t-1)}\mbox{[Eq.~(\ref{ap4:ddty})]},
\end{eqnarray*}
where the second term is added to eliminate second order partial derivatives,
the next term --- to eliminate $\partial F/\partial x$, 
and the last term --- to eliminate $\partial F/\partial y$.
\qed\\

\noindent
Kato himself \cite{Kato95} derived equation (\ref{ap4:singeq})  
via a conversion to Pfaffian systems. The specialization transformation
of Pfaffian systems is particularly elegant, but Kato's method requires conversion
of a matrix differential equation to an ordinary differential equation.  
The specialization problem for  partial differential equations can be considered
as a restriction problem for $D$-modules \cite{OaTa}. However, available implementations 
of algorithms for $D$-modules worked very slowly on our examples. Border bases \cite{Border}
rather than Gr\"obner bases look promisingly suitable for our strategy. Towards the end of Section 1 we mentioned the simplest computational view of looking for a linear relation between full derivatives (in the specialized variable) expressed in terms of a basis of partial derivatives.

\section{Application: Dihedral hypergeometric functions}
\label{sec:dihedral}

In the first case of Theorem \ref{th:f2cases}, we have a general function of the form
$\displaystyle\hpg21{a,b}{c}{\,t^2}$. The corresponding differential equation is a general pull-back transformation of Euler's equation (\ref{eq:euler}) with respect to a degree 2 covering ramified above two (of its three) singular points. The pull-backed equation has 4 singularities in general; it has 3 (or less)
singularities when a local exponent difference below a ramification point is equal to $\frac12$.
Theorem \ref{th:f2cases} says the general pull-back equation has solutions expressible as an $F_2(x,y)$ function with $x+y=2$. As we mentioned, the $F_2$ series does not converge unless it terminates.

An  interesting case is when the parameters $b_1,b_2$ in the first case of Theorem \ref{th:f2cases}
are zero or negative integers, say $b_1=-k$, $b_2=-\ell$. 
The immediate conclusion is that the functions
\begin{equation} \label{eq:dih}
\hpg{2}{1}{\frac{a}{2},\,\frac{a+1}{2}+\ell\,}{\frac{1}{2}-k}{\,z\,} \quad\mbox{and}\quad
(1+\sqrt{z})^{-a}\app2{a; -k,-\ell}{-2k,-2\ell}{\frac{2\sqrt{z}}{1+\sqrt{z}},\frac2{1+\sqrt{z}}}
\end{equation}
and
satisfy the same second order ordinary equation, with respect to $z$. The monodromy group of the 
$\hpgo21(z)$ function is a dihedral group; the local exponent differences at $z=0$ and $z=\infty$ are  the half-integers $k+\frac12$, $\ell+\frac12$, respectively.  Functions with a dihedral monodromy group can be expressed as elementary functions, since after a quadratic pull-back transformation (ramified above $z=0$, $z=\infty$, in our case) the monodromy group of the pull-backed equation is a cyclic group, and hypergeometric solutions become simple power or logarithmic functions.  
An example of such an elementary expression is
\begin{eqnarray} \label{dihedr2}
\hpg{2}{1}{\frac{a}{2},\,\frac{a+1}{2}\,}{\frac{1}{2}}{\,z\,} & =
& \frac{(1-\sqrt{z})^{-a}+(1+\sqrt{z})^{-a}}{2}.
\end{eqnarray}
Relatedly,
\begin{eqnarray} \label{dihedr3}
\hpg{2}{1}{\!\frac{a+1}{2},\,\frac{a+2}{2}}{\frac{3}{2}}{\,z\,} &
= & \frac{(1-\sqrt{z})^{-a}-(1+\sqrt{z})^{-a}}{2\,a\,\sqrt{z}} \quad (a\neq 0),\\
\hpg{2}{1}{\,\frac{1}{2},\;1\,}{\frac{3}{2}}{\,z\,} &
= & \frac{\log(1+\sqrt{z})-\log(1-\sqrt{z})}{2\,\sqrt{z}}, \\ \label{dihedr1}
\hpg{2}{1}{\frac{a}{2},\,\frac{a+1}{2}\,}{a+1}{\,z\,} & = &
\left(\frac{1+\sqrt{1-z}}{2} \right)^{-a}.
\end{eqnarray}
The $\hpgo21$ function in (\ref{eq:dih}) is a contiguous family of Gauss hypergeometric equations
with the same monodromy group. 

The relation between the two functions in (\ref{eq:dih}) allows us to generalize formula (\ref{dihedr2}).
The generalization is
\begin{eqnarray} \label{eq:dih12}
\!\!\frac{\left(\frac{a+1}2\right)_\ell}{\left(\frac12\right)_\ell}\,
\hpg{2}{1}{\frac{a}{2},\frac{a+1}{2}+\ell}{\frac{1}{2}-k}{\,z\,}\equal\frac{(1+\sqrt{z})^{-a}}2
\app2{a; -k,-\ell}{-2k,-2\ell}{\frac{2\sqrt{z}}{1+\sqrt{z}},\frac2{1+\sqrt{z}}}\nonumber\\
&\!\!\!+\!\!\!&\frac{(1-\sqrt{z})^{-a}}2
\app2{a; -k,-\ell}{-2k,-2\ell}{\frac{2\sqrt{z}}{\sqrt{z}-1},\frac2{1-\sqrt{z}}}.
\end{eqnarray}
The $F_2$ series on the right-hand side are 
finite sums with $(k+1)(\ell+1)$ terms.  If \mbox{$\left(\frac{a+1}2\right)_\ell\neq 0$},
this is an expression of the $\hpgo21$ function as a linear combination of two explicit solutions of
the same differential equation. The identity can be proved by using the symmetry
with respect to the conjugation of $\sqrt{z}$ and checking the value of both sides at $z=0$,
which leads to evaluation of $\displaystyle\hpg21{a,-\ell\,}{-2\ell}{\,2}$ obtainable by Zeilberger's algorithm. See \cite{Vidunas08} for details.
  
Similarly, a generalization of (\ref{dihedr3}) is
\begin{eqnarray} \label{eq:dih32}
\frac{\left(\frac{a+1}2\right)_k\left(\frac{a}2\right)_{k+\ell+1}}
{\left(\frac12\right)_k\left(\frac12\right)_{k+1}\left(\frac12\right)_\ell}\,(-1)^k\,z^{k+\frac12}\,
\hpg{2}{1}{\frac{a+1}{2}+k,\,\frac{a}{2}+k+\ell+1\,}{\frac{3}{2}+k}{\,z\,}\hspace{-160pt}\nonumber\\
\equal \frac{(1+\sqrt{z})^{-a}}2
\app2{a; -k,-\ell}{-2k,-2\ell}{\frac{2\sqrt{z}}{1+\sqrt{z}},\frac2{1+\sqrt{z}}}\nonumber\\
&\!\!\!-\!\!\!&\frac{(1-\sqrt{z})^{-a}}2
\app2{a; -k,-\ell}{-2k,-2\ell}{\frac{2\sqrt{z}}{\sqrt{z}-1},\frac2{1-\sqrt{z}}}.
\end{eqnarray}
To show this identity, one can verify that both sides satisfy the same recurrence relations with
respect to $k$ and $\ell$, and check the identity for a few values $(k,\ell)$.

Recalling the relation between $F_3(x,x/(2x-1))$ and $\hpgo21$ functions in Theorem \ref{th:f3cases}, 
Karlsson's identity \cite[9.4.(90)]{Srivastava85} gives the following generalization of (\ref{dihedr1}):
\begin{eqnarray} \label{eq:diha2}
\!\! \hpg{2}{1}{\frac{a-\ell}{2},\,\frac{a+\ell+1}{2}}{a+k+1}{\,z\,}\equal
\left(\frac{1+\sqrt{1\!-\!z}}{2} \right)^{-a-k}(1-z)^{k/2} 
\times\nonumber\\
&&\app3{k+1,\ell+1; -k,-\ell}{a+k+1}{\frac{\sqrt{1\!-\!z}-1}{2\sqrt{1-z}},\frac{1-\sqrt{1\!-\!z}}2}.
\end{eqnarray}
Note that the $F_3$ sum is finite for any integers $k,\ell$, since it is invariant under the substitutions
$k\mapsto-k-1$ and $\ell\mapsto-\ell-1$. 
The $F_2$ sums in (\ref{eq:dih12}) and (\ref{eq:dih32}) can be written as $F_3$ sums, or vice versa in (\ref{eq:diha2}), by reversing the double summation in both directions:
\begin{eqnarray}
\app2{a; -k,-\ell}{-2k,-2\ell}{x,y}=\frac{k!\,\ell!\,(a)_{k+\ell}}{(2k)!\,(2l)!}\,x^k\,y^{\ell}\,
\app3{k+1,\ell+1; -k,-\ell}{1-a-k-\ell}{\frac1x,\frac1y}.
\end{eqnarray}
These explicit expressions are analyzed more thoroughly in \cite{Vidunas08}.

\section{Appendix}

Here we recall a couple of relevant facts about ordinary Fuchsian equations, 
and review related explicit results on Appell' functions existing in literature.

The third order hypergeometric differential equation for 
$\hpg32{A,B,C}{D,\,E}{x}$ can be written 
\begin{eqnarray} \label{eq:hpgde3}
x^2(1-x)\,\frac{d^3y(x)}{dx^3}+
x\big(D+E+1-(A+B+C+3)x\big)\frac{d^2y(x)}{dx^2}+\nonumber \\
\big(DE-(AB\!+\!AC\!+\!BC\!+\!A\!+\!B\!+\!C\!+\!1)x\big)\frac{dy(x)}{dx}-A\,B\,C\,y(x)\equal 0.
\end{eqnarray}
The singularities and local exponents are:
\begin{eqnarray} \label{eq:hpg3le}
\begin{array}{llcl}
\mbox{at $x=0$}: & 0, & 1-D, & 1-E;\\
\mbox{at $x=1$}: & 0, & 1, &  D+E-A-B-C;\\
\mbox{at $x=\infty$}: & A, & B, & C.
\end{array}
\end{eqnarray}
Recall that a third order Fuchsian equation with 3 singularities is determined by the local exponents
and one {\em accessory parameter}. Say, one can add a scalar multiple of $y(x)/(x-1)$ to the left-hand side of equation (\ref{eq:hpgde3}) without changing the local exponents in (\ref{eq:hpg3le}).

Solutions of differential systems \mbox{(\ref{app2a})--(\ref{app2b})},
\mbox{(\ref{app3a})--(\ref{app3b})}, \mbox{(\ref{app4a})--(\ref{app4b})}, 
\mbox{(\ref{app1a})--(\ref{app1b})} in terms of bivariate hypergeometric series were considered by
many authors \mbox{\cite[pg. 222--242]{Erdelyi81}}, \mbox{\cite[pg.~291--305]{Srivastava85}}, 
starting from Le Vavasseur,  Appell, Kamp\'e de F\'eriet, 
Borng\"{a}sser, Erd\'elyi, Olsson. 

In general, each Appell's $F_2(x,y)$ function can be represented
by four $F_2$ series at the origin $(0,0)$:
\begin{eqnarray} \label{eq:f2s1}
\!\!\app2{a;\,b_1,b_2}{c_1,c_2}{x,y}
\equal  (1-x)^{-a}\app2{a;\,c_1-b_1,b_2}{c_1,c_2}{\frac{x}{x-1},\frac{y}{1-x}}\\
\equal  (1-y)^{-a}\app2{a;\,b_1,c_2-b_2}{c_1,c_2}{\frac{x}{1-y},\frac{y}{y-1}}\\ \label{eq:f2s4}
\equal  (1\!-\!x\!-\!y)^{-a}\app2{a;c_1\!-\!b_1,c_2\!-\!b_2}{c_1,c_2}
{\frac{x}{x\!+\!y\!-\!1},\frac{y}{x\!+\!y\!-\!1}}.
\end{eqnarray} 
Note that the parameters of the $F_2(x,2-x)$ function in (\ref{eq:f2f0}) are invariant 
under these transformations; this invariance is discussed in \cite[Section 10]{SasakiYoshida}.
Besides, in general there are four distinct $F_2(x,y)$ functions that are local solutions 
of \mbox{(\ref{app2a})--(\ref{app2b})} at the origin. The other 3 functions are represented by
the series
\begin{eqnarray} \label{eq:f2s5}
x^{1-c_1}\app2{1+a-c_1;\,1+b_1-c_1,b_2}{2-c_1,c_2}{x,y},\\  
y^{1-c_2}\app2{1+a-c_2;\,b_1,1+b_2-c_2}{c_1,2-c_2}{x,y},\\  \label{eq:f2s7}
x^{1-c_1}y^{1-c_2}\app2{2+a-c_1-c_2;\,1+b_1-c_1,1+b_2-c_2}{2-c_1,2-c_2}{x,y}.
\end{eqnarray}
The $F_2$ system \mbox{(\ref{app2a})--(\ref{app2b})} has Horn's $H_2$ series local solutions at
\mbox{$(x,y)=(0,\infty)$} and \mbox{$(x,y)=(\infty,0)$}, and four distinct $F_3$ local solutions 
at $(x,y)=(\infty,\infty)$. The latter are obtainable by applying relation (\ref{app2app3}) 
to the four series in  (\ref{eq:f2s1})--(\ref{eq:f2s4}). 
Application of relation (\ref{app2app3}) to (\ref{eq:f2s5})--(\ref{eq:f2s7})
does not give different $F_3$ series, but only shows invariance of (\ref{appf3}) 
under the permutations $a_1\leftrightarrow b_1$ and $a_2\leftrightarrow b_2$ of upper parameters. 
Hypergeometric solutions of the $F_3$ system (\ref{app3a})--(\ref{app3b}) are described similarly.

The system (\ref{app4a})--(\ref{app4b}) has four generally different 
$F_4$ local solutions at $(x,y)=(0,0)$. Besides (\ref{appf4}) we have
\begin{eqnarray} \label{eq:a4sola}
x^{1-c_1}\app4{1+a-c_1;1+b-c_1}{2-c_1,c_2}{x,y},\quad
y^{1-c_2}\app4{1+a-c_2;1+b-c_2}{c_1,2-c_2}{x,y},\\
x^{1-c_1}y^{1-c_2}\app4{2+a-c_1-c_2;2+b-c_1-c_2}{2-c_1,2-c_2}{x,y}.
\end{eqnarray}
Similar sets of four local solutions exist at $(x,y)=(0,\infty)$ and $(\infty,0)$. 
A connection formula between (\ref{appf4}) and the following two functions
is presented in \cite[pg 26]{Appel26} and 
\cite[9.4.(69)]{Srivastava85}, for example:
\begin{equation} \label{eq:a4solz}
y^{-a}\app4{a;1+a-c_2}{c_1,1+a-b}{\frac{x}{y},\frac1y},\qquad
y^{-b}\app4{1+b-c_2,b}{c_1,1+b-a}{\frac{x}{y},\frac1y}.
\end{equation}

The full set of $F_1$ solutions to the system  \mbox{(\ref{app1a})--(\ref{app1b})} consists 
of 60 $F_1$ series in general. First of all, the $F_1$ solutions are identified in
sextets: 
\begin{eqnarray} \label{eq:f1s1}
\!\!\!\app1{a; b_1, b_2}{c}{x,y} 
\equal(1-x)^{-b_1}(1-y)^{-b_2}\app1{c-a; b_1, b_2}{c}{\frac{x}{x-1},\frac{y}{y-1}},\\
\equal (1-x)^{-a}\app1{a; c-b_1-b_2,b_2}{c}{\frac{x}{x-1},\frac{x-y}{x-1}},\\
\equal (1-y)^{-a}\app1{a; b_1,c-b_1-b_2}{c}{\frac{x-y}{1-y},\frac{y}{y-1}},\\
\equal (1\!-\!x)^{c-a-b_1}(1\!-\!y)^{-b_2}\app1{c\!-\!a; c\!-\!b_1\!-\!b_2,b_2}{c}{x,\frac{x\!-\!y}{1\!-\!y}},\\
\equal (1\!-\!x)^{-b_1}(1\!-\!y)^{c-a-b_2}\app1{c\!-\!a; b_1,c\!-\!b_1\!-\!b_2}{c}{\frac{x\!-\!y}{x\!-\!1},y}.
\end{eqnarray}
Consequently, there are 10 different $F_1$ solutions of \mbox{(\ref{app1a})--(\ref{app1b})} in general.
The other 9 are represented by the following series:
\begin{eqnarray}
\app1{a; b_1, b_2}{1\!+\!a\!+\!b_1\!+\!b_2\!-\!c}{1\!-\!x,1\!-\!y}, \quad
x^{-b_1}y^{-b_2}\app1{1+b_1+b_2-c;b_1,b_2}{1+b_1+b_2-a}{\frac1x,\frac1y},\\
x^{-a}\app1{a; 1+a-c,b_2}{1+a-b_1}{\frac1x,\frac{y}x},\hspace{60pt}
y^{-a}\app1{a; b_1,1+a-c}{1+a-b_2}{\frac{x}y,\frac1y},\\
(1-x)^{-b_1}(1-y)^{c-a-b_2}\app1{c-a; b_1, c-b_1-b_2}{c-a-b_2+1}{\frac{1-y}{1-x},1-y},\\
(1-x)^{c-a-b_1}(1-y)^{-b_2}\app1{c-a; c-b_1-b_2, b_2}{c-a-b_1+1}{1-x,\frac{1-x}{1-y}},\\
x^{-b_1}y^{b_1-c+1}\app1{1+b_1+b_2-c;b_1,1+a-c}{2+b_1-c}{\frac{y}x,y},\\
x^{b_2-c+1}y^{-b_2}\app1{1+b_1+b_2-c;1+a-c,b_2}{2+b_2-c}{x,\frac{x}y},\\ \label{eq:f1s9}
x^{b_1+b_2-c}(x\!-\!y)^{1-b_1-b_2}(1\!-\!x)^{c-a-1}
\app1{1\!-\!b_1;1\!+\!a\!-\!c,c\!-\!b_1\!-\!b_2}{2-b_1-b_2}{\frac{x\!-\!y}{x\!-\!1},\frac{x\!-\!y}x}.
\end{eqnarray}
Besides, Appel's $F_1$ functions can be expressed as $F_3$ series following formula (\ref{eq:f1f3}). 
Appell's $F_2$ series realize $F_1$ functions as well; the following identity hold in the neighborhood
of the point $(x/y,y)=(1,0)$ of the blow-up at $(x,y)=(0,0)$:
\begin{eqnarray}
\app1{a; b_1, b_2}{c}{x,y}=(y/x)^{b_1}\app2{b_1+b_2;b_1,a}{b_1+b_2,c}{\frac{x-y}{x},y}.
\end{eqnarray}

In the remainder of the Appendix, we summarize identities in \cite[Section 9.4]{Srivastava85}
that relate different Appell's functions to each other or to Gauss hypergeometric functions.
The expressions for Appell's functions with reducible (up to possibly a quadratic transformation)
monodromy groups are formulas \mbox{\cite[9.4, (87), (88), (98), (99), (108)--(114)]{Srivastava85}};
several of these formulas are attributed to Bailey or to \cite{Appel26}.
It looks convenient to substitute $x\mapsto x/(x-1)$, $y\mapsto y/(y-1)$ in 
Bailey's formulas \mbox{\cite[9.4, (110)--(113)]{Srivastava85}}. We mentioned special
relations (\ref{eq:f2f4rel})--(\ref{eq:f2f4rel2}) between $F_2$ and $F_4$ functions; 
another similar relation is \mbox{\cite[9.4, (216)]{Srivastava85}}. The following 5 formulas 
are \mbox{\cite[9.4, (190)--(192), (149), (179)]{Srivastava85}}; they
give a taste of what can be expected from an exhaustive research of
univariate specializations of Appell's functions to $\hpgo43$ or $\hpgo32$ functions. 
The first three formulas are attributed to Bailey, and (\ref{eq:burchnal}) is due to Burchnall:
\begin{eqnarray}
\app2{a;b,b}{c,c}{x,-x}\equal \hpg43{\frac{a}2,\frac{a+1}2,b,c-b}
{c,\frac{c}2,\frac{c+1}2}{x^2},\\
\app2{a;b_1,b_2}{2b_1,2b_2}{x,-x}\equal 
\hpg43{\frac{a}2,\frac{a+1}2,\frac{b_1+b_2}2,\frac{b_1+b_2+1}2}
{b_1+\frac12,b_2+\frac12,b_1+b_2}{x^2},\\
\app3{a,a;b,b}{c}{x,-x}\equal \hpg43{a,b,\frac{a+b}2,\frac{a+b+1}2}
{a+b,\frac{c}2,\frac{c+1}2}{x^2},\\ \label{eq:burchnal}
\app4{a;b}{c_1,c_2}{x,x}\equal \hpg43{a,b,\frac{c_1+c_2}2,\frac{c_1+c_2-1}2}
{c_1,c_2,c_1+c_2-1}{4x},\\
\app4{a;b}{c,c}{x,-x}\equal \hpg43{\frac{a}2,\frac{a+1}2,\frac{b}2,\frac{b+1}2}
{c,\frac{c}2,\frac{c+1}2}{-4x^2}.
\end{eqnarray}


\end{document}